\documentclass[11pt]{amsart}

\input{diagrams}

\newtheorem{proposition}{Proposition}[section]
\newtheorem{lemma}[proposition]{Lemma}
\newtheorem{corollary}[proposition]{Corollary}
\newtheorem{theorem}[proposition]{Theorem}

\theoremstyle{definition}
\newtheorem{definition}[proposition]{Definition}
\newtheorem{example}[proposition]{Example}
\newtheorem{examples}[proposition]{Examples}
\newtheorem{remark}[proposition]{Remark}
\newtheorem{remarks}[proposition]{Remarks}

\newcommand{\thlabel}[1]{\label{th:#1}}
\newcommand{\thref}[1]{Theorem~\ref{th:#1}}
\newcommand{\selabel}[1]{\label{se:#1}}
\newcommand{\seref}[1]{Section~\ref{se:#1}}
\newcommand{\lelabel}[1]{\label{le:#1}}
\newcommand{\leref}[1]{Lemma~\ref{le:#1}}
\newcommand{\prlabel}[1]{\label{pr:#1}}
\newcommand{\prref}[1]{Proposition~\ref{pr:#1}}
\newcommand{\colabel}[1]{\label{co:#1}}
\newcommand{\coref}[1]{Corollary~\ref{co:#1}}
\newcommand{\relabel}[1]{\label{re:#1}}

\newcommand{\exlabel}[1]{\label{ex:#1}}
\newcommand{\exref}[1]{Example~\ref{ex:#1}}
\newcommand{\delabel}[1]{\label{de:#1}}

\newcommand{\eqlabel}[1]{\label{eq:#1}}
\newcommand{\equref}[1]{(\ref{eq:#1})}

\newcommand{\Hom}{{\rm Hom}}
\newcommand{\End}{{\rm End}}

\def\Tr{{\rm Tr}}

\def\ot{\otimes}

\def\sq{\square}

\def\ZZ{{\mathbb Z}}

\newcommand{\Cc}{\mathcal{C} }
\newcommand{\Dd}{\mathcal{D} }
\newcommand{\Ee}{\mathcal{E}}
\newcommand{\Ff}{\mathcal{F}}
\newcommand{\Gg}{\mathcal{G}}
\newcommand{\Hh}{\mathcal{H}}
\newcommand{\Mm}{\mathcal{M}}
\newcommand{\Pp}{\mathcal{P}}

\newcommand{\Tt}{\mathcal{T}}
\newcommand{\Kk}{\mathbb{K}}
\def\*C{{}^*\hspace*{-1pt}{\Cc}}
\def\text#1{{\rm {\rm #1}}}

\def\dul#1{\underline{\underline{#1}}}
\def\Nat{\dul{\rm Nat}}

\begin{document}
\title[Naturally full functors in nature]{Naturally full functors in nature}
\author{A. Ardizzoni}
\address{Departament of Mathematics, University of Ferrara, Via
Machiavelli 35, I-44100 Ferrara, Italy}
\email{ardiz@dm.unife.it}
\author{S. Caenepeel}
\address{Faculty of Engineering Sciences,
Vrije Universiteit Brussel, VUB, B-1050 Brussels, Belgium}
\email{scaenepe@vub.ac.be}
\urladdr{http://homepages.vub.ac.be/\~{}scaenepe/}
\author{C. Menini}
\address{Departament of Mathematics, University of Ferrara, Via
Machiavelli 35, I-44100 Ferrara, Italy}
\email{men@dns.unife.it}
\author{G. Militaru}
\address{Faculty of Mathematics, University of Bucharest, Str.
Academiei 14, RO-70109 Bucharest 1, Romania}
\email{gmilit@al.math.unibuc.ro}
\thanks{S. Caenepeel and G. Militaru were supported by the bilateral project
``Hopf Algebras in Algebra, Topology, Geometry and Physics" of the
Flemish and Romanian governments. This paper was written while C.
Menini was member of G.N.S.A.G.A. with partial financial support
from M.I.U.R.. G. Militaru was partially supported by I.N.D.A.M.,
while he was visiting professor at University of Ferrara.}


\begin{abstract}
We introduce and discuss the notion of naturally full functor. The definition
is similar to the definition of separable functor: a naturally full functor is a
functorial version of a full functor, while a separable functor is a functorial
version of a faithful functor. We study general properties of naturally full
functors. We also discuss when functors between module categories and
between categories of comodules over a coring are naturally full.
\end{abstract}

\maketitle

\section*{Introduction}
Separable functors were introduced by N\u ast\u asescu, Van den Bergh and
Van Oystaeyen \cite{NastasescuVV89}. The terminology is inspired by the
result that the restriction of scalars functor associated to a morphism of rings
is a separable functor if and only if the corresponding ring extension is
separable. Separable functors have been studied extensively
during the passed decade, we refer to \cite{CMZ2002} for a detailed discussion
of results and applications. The original definition from \cite{NastasescuVV89}
can be restated in a more categorical way as follows: to a functor
$F:\ \Cc\to \Dd$, we can associate a natural transformation
$\Ff:\ \Hom_{\Cc}(\bullet,\bullet)\to \Hom_{\Dd}(F(\bullet),F(\bullet))$,
mapping $f$ to $F(f)$; $F$ is called separable if $\Ff$ has a left inverse.
Recall that $F$ is called faithful if every $\Ff_{C,C'}$ is injective, i.e.
it has a left inverse in the category of sets. Therefore a separable functor
is faithful; in fact we could call a separable functor naturally faithful, in the
sense that every $\Ff_{C,C'}$ has a left inverse that is functorial in $C$ and
$C'$.\\
The aim of this note is to study the notion of naturally full
functors; a functor $F$ is naturally full if $\Ff$ has a right
inverse; if $F$ is full, then every $\Ff_{C,C'}$ is surjective,
and has a right inverse - $F$ is
naturally full when this right inverse is functorial.\\
In \seref{2}, we study general properties of naturally full functors; some of
these properties are analogous to some properties of separable functors:
We have full versions of the Rafael and Rafael-Frobenius Theorems telling
when a functor having an adjoint (resp. a Frobenius functor) is naturally full.\\
In \seref{3}, we study some particular examples. From an algebraic
point of view, the first example to look at is extension and
restriction of scalars. The restriction of scalars functors is
naturally full if and only if it is full. The same result is not
true for the extension of scalars functor, see \exref{3.2}. We can
extend our results to module categories connected by a bimodule,
see \seref{3.2}. In Sections \ref{se:3.3} and \ref{se:3.4}, we look at
categories of comodules over a coring.

\section{Prelimaries}\selabel{1}
We use the following conventions. For an object $V$ in a category,
the identity morphism $V\to V$ is denoted by $V$. A functor is assumed to
be covariant. By a ring, we will always
mean a ring with unit.
For a ring $R$,  $\Mm _R$ (resp. ${}_R\Mm $, ${}_{R}\Mm _{R}$) denotes
the category of right $R$-modules (resp.\ left $R$-modules,
$R$-bimodules).

\subsection{Adjoint functors}\label{subs.adj}\selabel{1.1}
Let $(F, G)$ be a pair of adjoint functors between two categories
$\Cc$ and $\Dd$. This means that there exist two
natural transformations $\eta:\ 1_\Cc\to GF$
and $\varepsilon:\ FG\to 1_\Dd$, called the unit and
counit of the adjunction, such that
\begin{equation}\label{1.1}
G(\varepsilon_D)\circ\eta_{G(D)}= I_{G(D)}~~~~{\rm and}~~~~
\varepsilon_{F(C)}\circ F(\eta_C)= I_{F(C)}
\end{equation}
for all $C\in \Cc$ and $D\in \Dd$. Then there exists
a natural isomorphism
\begin{equation}\label{1.2a}
\varphi_{C, D}: \Hom_{\mathcal D} (FC, D)\to
\Hom_{\Cc }(C, GD), \quad
\varphi_{C, D} (f) = G(f) \circ \eta_C
\end{equation}
with inverse
\begin{equation}\label{1.2b}
\psi_{C, D} = \varphi_{C, D}^{-1}:
\Hom_{\Cc }(C, GD)\to
\Hom_{\Dd} (FC, D), \quad
\psi_{C, D} (g) =\varepsilon_D \circ F(g)
\end{equation}
for any $C\in \Cc$ and $D\in \Dd$. Conversely,
$\eta$ and $\varepsilon$ are recovered from
$\varphi$ and $\psi$ by
$$
\eta_C = \varphi_{C, FC} (FC), \quad
\varepsilon_D = \psi_{GD, D} (GD).
$$
We also have isomorphisms (see for example \cite[Proposition 11]{CMZ2002}):
$$
\alpha:\  \Nat(GF, 1_{\Cc })\to
\Nat(\Hom_{\Dd}(F(\bullet), F(\bullet)),
\Hom_{\Cc }(\bullet,\bullet))
$$
$$
\beta: \Nat(1_{\Dd}, FG)\to
\Nat(\Hom_{\Cc }(G(\bullet), G(\bullet)),
\Hom_{\Dd}(\bullet,\bullet))
$$
defined as follows:
\begin{equation}\label{1.41}
\alpha(\nu)_{C,C'}(g) = \nu_{C'}\circ G(g)\circ \eta_C
\end{equation}
for all $\nu \in \Nat(GF, 1_{\mathcal C})$, $C$,
$C' \in {\Cc }$ and $g:\ F(C)\to F(C')$ in
${\Dd}$ with  inverse
\begin{equation}\label{1.5}
\alpha^{-1}({\Pp})_{C} =
{\Pp}_{GFC, C} (\varepsilon_{FC})
\end{equation}
for any
${\Pp} \in \Nat(\Hom_{\mathcal D}(F(\bullet), F(\bullet)),
\Hom_{\mathcal C}(\bullet,\bullet))$ and $C\in {\mathcal C}$ and
respectively,
\begin{equation}\label{1.6}
\beta(\xi)_{D, D'}(f) = \varepsilon_{D'}\circ F(f)\circ \xi_D
\end{equation}
for any $\xi \in \Nat(1_{\mathcal D}, FG)$, $D$,
$D' \in {\mathcal D}$ and $f:\ G(D)\to G(D')$ in
${\mathcal C}$ with inverse
\begin{equation}\label{1.7}
\beta^{-1}({\mathcal T})_{D} = {\mathcal T}_{D, FGD} (\eta_{GD})
\end{equation}
for any
${\Tt }\in \Nat(\Hom_{\mathcal C}(G(\bullet), G(\bullet)),
\Hom_{\mathcal D}(\bullet,\bullet))$ and $D\in {\mathcal D}$.\\
Assume now that $G$ is a left adjoint of $F$ and let
$$\mu:\ GF\to 1_{\Cc}~~{\rm and}~~\chi:\ 1_{\Dd}\to FG$$
be the counit and unit of the adjunction $(G, F)$. By
\cite[Proposition 10]{CMZ2002},there exist one-to-one
corespondences between the following classes
$$\Nat(GF, 1_{\mathcal C})\cong \Nat(G, G)\cong \Nat(F, F)
\cong \Nat(1_{\mathcal D}, FG)
$$
Let us describe these correspondences: for any natural transformation $\nu  :\ GF \to
1_{\Cc }$ there exist unique natural transformations
$\alpha : \ G\to G$, $\beta:\ F\to F$, such that
\begin{equation}\label{1.3}
\mu_{C}\circ \alpha_{FC} = \nu_{C} = \mu_{C}\circ G(\beta_C);
\end{equation}
for any natural transformation $\xi :\ 1_{\mathcal D} \to FG$,
there exist unique natural transformations $\alpha :\ G\to G$,
$\beta:\ F\to F$, such that
\begin{equation}\label{1.4}
\beta_{GD}\circ \chi_{D} = \xi_{D} =
F(\alpha_D)\circ \xi_{D}.
\end{equation}

Recall that a functor $F$ is called Frobenius if it has a right
adjoint $G$ which is also a left adjoint. In this case, there
exist four functors $\varepsilon$, $\eta$, $\xi$ and $\mu$
satisfying the properties given above.

\subsection{Separable functors}\label{ffs}\selabel{1.2}
Separable functors were introduced in \cite {NastasescuVV89}; several applications
have appeared in the literature, see \cite{CMZ2002}. The definition of separable functor
can be formulated in the following way: for a functor
$F :\  \Cc \to \Dd$, we gave two functors
$$\Hom_{\Cc}(\bullet,\bullet),
~~~\Hom_{\Dd}(F(\bullet), F(\bullet)):\  \Cc^{\rm op}\times
\Cc\to\dul{\rm Sets}$$
and a natural transformation
$$
{\mathcal F}:\  \Hom_{\Cc}(\bullet, \bullet)\to
\Hom_{\Dd}(F(\bullet), F(\bullet)),
\quad {\mathcal F}_{C,C'}(f) =F(f)
$$
for any $f:\ C\to C'$ in $\Cc$. The natural transformation
${\mathcal F}$ will play a key role in the sequel.
$F$ is called {\it separable} if
${\mathcal F}$ splits, that is, we have a natural
transformation
$$
{\mathcal T}:\  \Hom_{\Dd}(F(\bullet), F(\bullet)) \to
\Hom_{\Cc}(\bullet, \bullet)
$$
such that
$${\mathcal T}\circ {\mathcal F} =
{\bf 1}_{\Hom_{\Cc}(\bullet, \bullet)}$$
If $F$ is separable, then for all $C$, $C'\in {\Cc}$, the map
${\mathcal F}_{C,C'}:\ \Hom_{\Cc} (C,C') \to \Hom_{\mathcal
D}(F(C), F(C'))$ is injective, since it has a left inverse, and it
follows that $F$ is a faithful functor. Since an injective map between
sets has a left inverse, we can restated the definitions of separable
and faithful functors in the following way: $F$ is faithful if every
${\mathcal F}_{C,C'}$ has a left inverse; $F$ is separable if this
left inverse can be chosen to be natural in $C$ and $C'$.
Perhaps it would be better, at least from the categorical point of
view, to call separable functors naturally faithful functors.
Observe that a faithful functor is not necessarily separable:
if $K/L$ is a purely inseparable field extension, then the restriction
of scalars functor is faithful, but not separable.\\

Now recall the dual notion of a faithful functor: $F$ is called a
full functor if every ${\mathcal F}_{C,C'}$ is surjective, or,
equivalently, has a right inverse. As far as we know, there are
not many characterizations of full functors avalaible in the
literature. In the case where $F$ is the restriction of scalars
functor, we have the following result (see for example
\cite[Proposition XI.1.1 and XI.1.2]{St}). For $M\in {}_S\Mm_S$,
let $M^S = \{m\in M~|~sm = ms, ~{\rm for~all~}s\in S \}$ be the
set of $S$-invariant elements of $M$.

\begin{theorem}\label{th2}\thlabel{1.1}
Let $\varphi : R \to S$ be a ring morphism, and $\varphi_{*}:\
{}_S\Mm\to {}_R\Mm$ the restriction of scalars functor.
The following
statements are equivalent:
\begin{enumerate}
\item The restriction of scalars functor
$\varphi_{*} : \ {}_S\Mm \to {}_R\Mm$ is full;
\item $\varphi :\  R \to S$ is an epimorphism of rings i.e. it is an
epimorphism in the category of unital rings;
\item $M^R = M^S$, for any $M\in {}_S\Mm_S$;
\item $1_S\ot_R 1_S$ is a separability idempotent for the extension
$S/R$, i.e. $s\ot_R 1_S = 1_S \ot_R s$, for any $s\in S$;
\item the map
$$
\varepsilon_S : S\ot_R S \to S, \quad
\varepsilon_S (s_1\ot_R s_2) = s_1 s_2
$$
is injective (hence an isomorphism in ${}_S\Mm_S$);
\item the counit of the adjunction $(S\ot_R \bullet, \;\; \varphi_{*})$
$$
\varepsilon_N : S\ot_R N \to N, \quad
\varepsilon_N (s\ot_R n) = sn
$$
is an isomorphism of left $S$-modules, for all $N\in {}_S\Mm$.
\end{enumerate}
\end{theorem}

The fact that the
restriction of scalars functor $\varphi_{*} : {}_S\Mm \to
{}_R\Mm$ is a full functor explicitly means that
$$
{}_R\Hom (M, N) = {}_S\Hom (M, N)
$$
for any $M$, $N \in {}_S\Mm$, i.e. any left $R$-linear
map between two left $S$-modules is also left $S$-linear.

\section{Naturally full functors}\selabel{2}
We keep the notation of \seref{1}. Following the philosophy of \seref{1},
we introduce the notion of naturally full functor.

\begin{definition}\delabel{2.1}
A functor $F$ is called a naturally full functor,  if $\Ff$ cosplits. This means that there exists
a natural transformation
$$
{\Pp}:\  \Hom_{\Dd}(F(\bullet), F(\bullet)) \to
\Hom_{\Cc}(\bullet, \bullet)
$$
such that
\begin{equation}\eqlabel{2.1.1}
{\mathcal F}\circ {\Pp} =
{\bf 1}_{\Hom_{\Dd}(F(\bullet), F(\bullet))}
\end{equation}
\end{definition}

\begin{remarks}\relabel{2.2}
1. The fact that ${\Pp}$ is a natural transformation
means the following:  for any $X,Y,Z,T\in \Cc$ and
$f:\ X\to Y$, $h:\ Z\to T$ in $\Cc$ and
$g:\ F(Y)\to F(Z)$ in $\Dd$,
we have
\begin{equation}\eqlabel{2.2.1}
{\Pp}_{X,T} (F(h)\circ g\circ F(f))=
h\circ {\Pp}_{Y,Z}(g)\circ f.
\end{equation}
\equref{2.1.1} can be rewritten as
\begin{equation}\eqlabel{2.2.2}
F( {\Pp}_{C,C'}(u) ) = u
\end{equation}
for all $C,C' \in \Cc$ and $u:\ FC\to FC'$ in $\Dd$.\\
2. A naturally full functor is always full, but the converse is not true in
general. For a counterexample, see \exref{3.2}.\\
3. The following conditions are equivalent:
\begin{enumerate}
\item $F$ is fully faithful; \item $F$ is separable and naturally
full.
\end{enumerate}
4. Combining 3. with the well-known result that a functor is an equivalence if and
only if it is fully faithful and surjective on objects, we find that
the following conditions are equivalent:
\begin{enumerate}
\item $F$ is a category equivalence; \item $F$ is fully faithful
and surjective on objects; \item $F$ is separable, naturally full
and surjective on objects.
\end{enumerate}
\end{remarks}

We will now present some general properties of naturally functors. The proofs are
sometimes similar to corresponding proofs for properties of separable functors,
see \cite{CMZ2002}.

\begin{proposition}\prlabel{2.3}
Consider functors $G:\ \Cc \to \Dd$ and $H:\ \Dd \to \Ee$.
\begin{enumerate}
\item If $G$ and $H$ are naturally full, then $H\circ G$ is also naturally full.
\item If $H\circ G$ is a naturally full and $H$ is faithful, then
$G$ is naturally full.
\end{enumerate}
\end{proposition}

\begin{proof}
1) Obvious: if ${\Pp^G}$ and ${\Pp^H}$ are right inverses of $\Gg$ and $\Hh$,
then $\Pp^G\circ \Pp^H$ is a right inverse of $\Hh\circ \Gg$.\\
2) Let ${\Pp^{H\circ G}}$ be a right inverse of $\Hh\circ \Gg$.
Then
$$\Hh\circ \Gg\circ \Pp^{H\circ G}\circ \Hh=\Hh,$$
and from the fact that $H$ is faithful, it follows that $\Pp^{H\circ G}\circ \Hh$
is a right inverse of $\Gg$, as needed.
\end{proof}

\begin{proposition}\prlabel{2.4}
Consider functors $F:\ \Dd \to \Cc$, $G:\ \Cc \to \Dd$ and $H:\ \Dd \to
\Ee$ and assume that $G$ is a right or left adjoint of $F$. If $H\circ G$
is naturally full and $F$ is separable, then $H$ is naturally full.
\end{proposition}

\begin{proof}
Let $\Pp^{H\circ G}$ be a right inverse of $\Hh\circ \Gg$. If $G$ is a right
adjoint of $F$, then by Rafael's Theorem \cite{Rafael90} the unit
$\eta:\ 1_\Cc\to
GF$ of the adjunction splits, i.e. there exists a natural transformation $\nu :\ GF
\to 1_\Cc$ such that
$\nu_C\circ \eta_C  = C$,
for all $C \in \Cc$. We now define a natural transformation
$\mathcal P^{H}$ using (\ref{1.41}):
\begin{eqnarray*}
&&\hspace*{-2cm}
\mathcal P^{H}_{C,C'} (f)=\alpha(\nu)_{C,C'}\{{{\mathcal
P^{H\circ G}}[H(\eta_{C'})\circ f\circ
H(\nu_{C})]}\}\\
&=&\nu_{C'}\circ G\{{\mathcal P^{H\circ
G}}[H(\eta_{C'})\circ f\circ H(\nu_{C})]\}\circ \eta_{C}.
\end{eqnarray*}
for all $f \in {\Hom_{\Ee}(H(C), H(C'))}$. We easily compute that
\begin{eqnarray*}
&&\hspace*{-2cm}
{\Hh_{C,C'}}\circ {\mathcal P^{H}_{C,C'}}(f) =
H(\nu_{C'})\circ HG\{{\mathcal P^{H\circ
G}}[H(\eta_{C'})\circ f\circ H(\nu_{C})]\}\circ H(\eta_{C})\\
&=& H(\nu_{C'})\circ H(\eta_{C'})\circ f\circ H(\nu_{C})\circ
H(\eta_{C})= f,
\end{eqnarray*}
as needed.\\
The proof is similar in the situation where $G$ is a left adjoint
of $F$: by Rafael's Theorem, the counit $\varepsilon:\ GF\to
1_{\Dd}$ has a right adjoint $\xi$: $\varepsilon_D\circ \xi_D=D$,
for all $D\in\Dd$. We then define a natural transformation $\Pp^H$
using (\ref{1.6}).
\end{proof}

Rafael's Theorem \cite{Rafael90} provides an easy characterization of the separability of
a functor having an adjoint. Similar characterizations can be formulated for full and
faithful functor (\prref{2.5}) and for naturally full functors (\thref{2.6}).

\begin{proposition}\prlabel{2.5}
Let $G:\ \Dd\to \Cc$ be a right adjoint of $F:\ \Cc\to\Dd$.
\begin{enumerate}
\item $F$ is faithful if and only if $\eta_C$ is a monomorphism in $\Cc$, for all $C\in \Cc$;
\item $F$ is full if and only if $\eta_C$ cosplits in $\Cc$, for all $C\in \Cc$, that is, there
exists $\nu_C\in\Cc$ such that $\nu_C\circ\eta_C=GFC$;
\item $G$ is faithful if and only if $\varepsilon_D$ is a epimorphism in $\Dd$, for all $D\in \Dd$;
\item $G$ is full if and only if $\varepsilon_D$ splits in $\Dd$, for all $D\in \Dd$, that is, there
exists $\xi_D\in\Dd$ such that $\xi_D\circ\varepsilon_D=FGD$.
\end{enumerate}
\end{proposition}

\begin{proof}
For any $C,C'\in \Cc$, we consider the composition
$$\Omega_{C,C'}=\varphi_{C,FC'}\circ \Ff_{C,C'}:\
\Hom_{\Cc}(C,C')\to \Hom_{\Cc}(C,GFC'),$$
where $\varphi$ is the natural isomorphism defined in \seref{1.1}. We easily compute that
$$\Omega_{C,C'}(f)=\varphi_{C,C'}(F(f))=GF(f)\circ\eta_C=\eta_{C'}\circ f.$$
Since $\varphi$ is an isomorphism, injectivity or surjectivity of $\Ff_{C,C'}$ is
equivalent to injectivity or surjectivity of $\Omega_{C,C'}$.\\
Now $\Omega_{C,C'}$ is injective if and only if $\eta_{C'}$ is a monomorphism, and
this proves the first statement.\\
If $\Omega_{GFC,C}:\ \Hom_{\Cc}(GFC,C)\to \Hom_{\Cc}(GFC,GFC)$
 is surjective, then there exists a $\nu_C:\ GFC\to C$ such that
 $$\Omega_{GFC,C}(\nu_C)=\eta_C\circ\nu_C=GFC$$
 Conversely, assume that for every $C\in \Cc$, there exists $\nu_C$ such that
 $\eta_C\circ\nu_C=GFC$. For every $f:\ C\to GFC'$ in $\Cc$, we have
 $$\Omega_{C,C'}(\nu_{C'}\circ f)= \eta_{C'}\circ \nu_{C'}\circ f=f$$
 and it follows that $\Omega_{C,C'}$ is surjective. This finishes the proof of
 the second statement.\\
 The two remaining properties are the duals of the first and the second one.
\end{proof}

Our next result is Rafael's Theorem for naturally full functors.

\begin{theorem}\thlabel{2.6}
Let $(F,G)$ be an adjoint pair of functors between $\Cc$ and $\Dd$,
with unit $\eta$ and counit $\varepsilon$.\\
1. $F$ is naturally full if and only if $\eta:\ 1_\Cc\to GF$
cosplits, i.e. there exists a natural transformation $\nu :\ GF \to
1_\Cc$ such that
\begin{equation}\eqlabel{2.6.1}
\eta_C\circ \nu_C = {GFC}
\end{equation}
for all $C \in \Cc$. In this case, $\varepsilon_{FC}$ is an isomorphism in $\Dd$
with inverse $F(\eta_C)$.\\
2. $G$ is an naturally full if and only if $\varepsilon:\ FG\to
1_\Dd$ splits, i.e. there exists a natural transformation $\xi:\ 1_\Dd \to FG$ such that
\begin{equation}\label{2.6.2}
\xi_D\circ \varepsilon_D = {FGD}
\end{equation}
for all $D \in \Dd$. In this case, $G(\varepsilon_{D})$ is an isomorphism in $\Cc$
with inverse $\eta_{GD}$.
\end{theorem}

\begin{proof}
1. Assume first that $F$ is naturally full and let $\Pp$ be a right inverse of $\Ff$.
 Let
$\nu= \alpha^{-1}({\mathcal P})$ be the natural transformation
given by (\ref{1.5}):
$$
\nu_C : \ GFC \to C, \quad
\nu_C = {\mathcal P}_{GFC, C} (\varepsilon_{FC})
$$
We consider the natural isomorphism  $\psi$ given by (\ref{1.2b}).
Then we have
\begin{eqnarray*}
 &&\hspace*{-15mm}
\psi (\eta_C \circ \nu_C)\stackrel{(\ref{1.2b}) } {=}
\varepsilon_{FC}\circ F(\eta_C \circ \nu_C)
= \varepsilon_{FC}\circ F(\eta_C) \circ
F({\mathcal P}_{GFC, C} (\varepsilon_{FC}) )\\
&\stackrel{\equref{2.2.2} } {=}&
\varepsilon_{FC}\circ F(\eta_C) \circ \varepsilon_{FC}
\stackrel{(\ref{1.1}) } {=}
\varepsilon_{FC} = \varepsilon_{FC}\circ F(Id_{GFC})
\stackrel{(\ref{1.2b}) } {=}
\psi (Id_{GFC})
\end{eqnarray*}
From the fact that  $\psi$ is a natural isomorphism, we obtain that
$\eta_C \circ \nu_C = Id_{GFC}$, for all $C \in {\mathcal C}$.
Furthermore,
$$
{GFC} =  \eta_C \circ \nu_C = \eta_C \circ {\mathcal P}_{GFC, C}
(\varepsilon_{FC}) \stackrel{\equref{2.2.1} } {=} {\mathcal
P}_{GFC,GFC} (F(\eta_C)\circ \varepsilon_{FC}).$$ Applying $F$ and
using \equref{2.2.2}, we obtain that $F(\eta_C)\circ
\varepsilon_{FC} = Id_{FGFC}$. Combining this with (\ref{1.1}), we
see that
$\varepsilon_{FC}$ and $F(\eta_C)$ are each others inverses.\\
Conversely, let $\nu :\ GF \to 1_\Cc$ be such that $\eta_C\circ \nu_C
= {GFC}$ for all $C \in \Cc$ and let ${\mathcal P}
=\alpha (\nu)$ be the natural transformation given by
$(\ref{1.41})$, i.e.
$$
{\mathcal P}_{C, C'} (f) = \nu_{C'} \circ G(f) \circ \eta_C
$$
for all $C$, $C' \in {\mathcal C}$ and $f:\ F(C)\to F(C')$ in
${\mathcal D}$. We have to show that
$F( {\mathcal P}_{C,C'}(f) ) = f$. We easily compute that
\begin{eqnarray*}
&&\hspace*{-2cm}
F(\nu_{C'} \circ G(f) \circ \eta_C) =
Id_{FC'} \circ F(\nu_{C'} \circ G(f) \circ \eta_C)\\
&\stackrel{(\ref{1.1}) } {=}&
\varepsilon_{FC'} \circ F(\eta_{C'}) \circ
F(\nu_{C'} \circ G(f) \circ \eta_C)\\
&=&\varepsilon_{FC'} \circ
F(\eta_{C'}\circ \nu_{C'} \circ G(f) \circ \eta_C)\\
&=& \varepsilon_{FC'} \circ
F(G(f) \circ \eta_C )
\stackrel{(\ref{1.2a}) } {=}
(\psi \circ \varphi) (f) = f,
\end{eqnarray*}
as needed.\\
2. The second statement follows from the first one by duality arguments.
\end{proof}

The Rafael criterion for separability simplifies further if we consider a Frobenius
functor. In \prref{2.7}, we give necessary and sufficient conditions for the
natural fullness of a Frobenius functor.

\begin{proposition}\prlabel{2.7}
Let $F:\ \Cc\to \Dd$ be a Frobenius functor, with left and right adjoint $G$.
We use the same notation as in \seref{1} for the unit and counit of the two
adjunctions. Then the following statements are equivalent.
\begin{enumerate}
\item $F$ is a naturally full;
\item there exists a natural transformation $\alpha :\ G \to G$
such that
\begin{equation}\eqlabel{2.7.1}
\eta_C \circ \mu_C \circ \alpha_{FC} = {GFC},
\end{equation}
for all $C\in \Cc$;
\item there exists a natural transformation $\beta :\ F\to F$
such that
\begin{equation}\eqlabel{2.7.2}
\eta_C \circ \mu_C \circ G(\beta_C) ={GFC},
\end{equation}
for all $C\in \Cc$;
\item there exists a natural transformation $\tilde{\alpha} :\  G \to G$
such that
\begin{equation}\eqlabel{2.7.3}
\tilde{\alpha}_{FC} \circ \eta_C \circ \mu_C = {GFC},
\end{equation}
for all $C\in \Cc$;
\item there exists $\tilde{\beta} :\  F \to F$ a natural transformation
such that
\begin{equation}\eqlabel{2.7.4}
G(\tilde{\beta}_{C}) \circ \eta_C \circ \mu_C = {GFC},
\end{equation}
for all $C\in \Cc$.
\end{enumerate}
In this case, $\eta_C \circ \mu_C$ is an isomorphism in $\Cc$ and
$$
\alpha_{FC} = \tilde{\alpha}_{FC} = G(\beta_C) = G(\tilde{\beta}_{C}),
$$
\end{proposition}

\begin{proof}
1. Apply first \thref{2.6} to the adjunction $(F, G)$ and
then (\ref{1.3})to the adjunction $(G, F)$, in order to describe all
natural transformations $\nu$ in terms of $\alpha$ and $\beta$.
This entails the equivalence of 1), 2) and 3).\\
Then we view $F$ as a right adjoint of $G$ and apply \thref{2.6} and (\ref{1.4})
to the adjunction $(G, F)$; we obtain a description of all natural transformations
$\xi$ in terms of $\tilde{\alpha}$ and
$\tilde{\beta}$, and we find the equivalence 1), 4) and 5).
\end{proof}

\section{Applications and examples}\selabel{3}
\subsection{Extension and restriction of scalars}\selabel{3.1}
Let $\varphi:\ R\to S$ be a morphism of rings. The restriction of scalars functor
$\varphi_*:\ {}_S\Mm\to {}_R\Mm$ is a right adjoint of the extension of scalars functor
$\varphi^*=S\ot_R\bullet:\ {}_R\Mm\to {}_S\Mm$. The unit and counit of the adjunction
are  $\eta_M=\varphi\ot_R M:\ M\to S\ot_RM$ and $\varepsilon_N:\ S\ot_R N\to N$,
$\varepsilon_N(s\ot_R n)=sn$. It is well-known and easy to prove, see for
example \cite{CMZ2002}, that
$$\Nat(\varphi_*\circ\varphi^*,1_{{}_R\Mm})\simeq {}_R\Hom(S,R)_R~~{\rm and}~~
\Nat(1_{{}_S\Mm},\varphi^*\circ\varphi_*)\simeq (S\ot_RS)^S.$$
The natural transformation $\nu^E$ corresponding to $E\in {}_R\Hom(S,R)_R$
is given by
$$\nu^E_M:\ S\ot_RM\to M,~~\nu^E_M(s\ot_Rm)=E(s)m.$$
The natural transformation $\xi^e$ corresponding to $e=\sum e^1\ot e^2\in
(S\ot_RS)^S$ is given by
\begin{equation}\eqlabel{3.0.1}
\xi^e_N:\ N\to S\ot_RN,~~\xi^e_N(n)=\sum e^1\ot_Re^2n.
\end{equation}

\begin{proposition}\prlabel{3.1}
Let $\varphi : R\to S$ be a ring morphism.\\
1. The restriction of scalars functor $\varphi_{*}$
is naturally full if and only if it is full.\\
2. The following statements are equivalent:
\begin{enumerate}
\item The extension of scalars functor $\varphi ^{* }$ is naturally full;
\item there exists an $E\in {}_R\Hom (S, R){}_R$ such that $
\varphi\circ E={S}$;
 \item there exists a central idempotent $e$ of $R$ such that $S\cong Re$
 and $\varphi:\ R\to S\cong Re$ is the projection $\varphi(r)=re$.
\end{enumerate}
\end{proposition}

\begin{proof}
1. If $\varphi_*$ is full, then it follows from the equivalence $(1)\Longleftrightarrow (6)$
in \thref{1.1} that the counit $\varepsilon$ is a natural isomorphism; in particular
$\varepsilon$ has a left inverse, and it follows from \thref{2.6} that
$\varphi_*$ is naturally full.\\
2. $(1)\Longrightarrow(2)$. By \thref{2.6}, $\eta$ has a right
inverse $\nu$; take the corresponding $E\in {}_R\Hom(S,R)_R$. For
all $s\in S$, we have
$$s=(\eta_R\circ\nu^E_M)(s)=\eta_R(E(s))=\varphi(E(s))$$
$(2)\Longrightarrow(1)$. For all $M\in {_{R}\Mm}$, $m\in M$ and
$s\in S$, we compute
\begin{eqnarray*}
&&\hspace*{-2cm}
(\eta _{M}\circ \nu _{M}^{E})(s\otimes _{R}m)=\eta
_{M}(E(s)m)=1_{S}\otimes _{R}E(s)m\\
&=&1_{S}E(s)\otimes _{R}m=\varphi
(E(s))\otimes _{R}m=s\ot_R m.
\end{eqnarray*}
$(2)\Longrightarrow(3)$.
Let $e=E(1_S)$. Then
\begin{eqnarray*}
&&\hspace*{-2cm}
r\cdot e=rE(1_S)=E(r1_S)=E(\varphi(r)1_S)=E(\varphi(r))\\
&=&E(1_S\varphi)=E(1_Sr)=E(1_S)r=e\cdot r
\end{eqnarray*}
and
\begin{eqnarray*}
&&\hspace*{-2cm}
e^2=E(1_S)\cdot E(1_S)=E(E(1_S)\cdot
1_S)\\
&=&E(\varphi(E(1_S))1_S)=E(1_S1_S)=E(1_S)=e,
\end{eqnarray*}
so $e$ is a central idempotent in $R$. The restriction
$\varphi'=\varphi_{|Re}:\ Re\to S$ is an isomorphism, since
$$E(\varphi(re))=E(\varphi(r)\varphi(e))=rE(\varphi(e))=rE(1_S)=re$$
and $\varphi\circ E)=S$.\\
$(3)\Longrightarrow(2)$ is trivial.
\end{proof}

\begin{remarks}\relabel{3.1a}
1) It follows from \prref{3.1} that natural fullness of the restriction and extension of
scalars functor is left-right symmetric: $\varphi_*$ is naturally full if and only
if ${}_*\varphi:\ \Mm_S\to \Mm_R$ is naturally full.\\
2) Recall that a map of left $R$-modules $f:\ M\to M'$ is called pure
if the map
$$N\ot_R f:\ N\ot_R M\to N\ot_R M'$$
is injective, for every $N\in \Mm_R$. If $\varphi:\ R\to S$ is pure as a morphism
of left or right $R$-modules, then the map
$$R\to \{s\in S~|~s\ot_R 1=1\ot_R s\}$$
is an isomorphism (see for example \cite[Prop. 2.1]{Caenepeel03}).
If $S$ is faithfully flat as a left (right) $R$-module, then
$\varphi$ is pure as morphism
of left (right) $R$-modules.\\
Assume that $\varphi$ is left (or right) pure, and that $\varphi_*$ is naturally
full. Then it follows from condition (4) in \thref{1.1} that $\varphi$ is the
identity, hence $R=S$.\\
3) Assume that $\varphi^*$ is naturally full. It follows from \prref{3.1} that,
up to isomorphism, we may assume that $S=Re$, where $e$ is a central idempotent
of $R$. Then $S\ot_R S=Re\ot_R Re\cong Re$, and $e\ot e$ is a separability idempotent
of $S/R$. Now $e$ is the unit element of $S$, so $\varphi_*$ is also naturally full.
Furthermore, the inclusion $i:\ S\to R$ is an $R$-bimodule map, and
$(e\ot_R e,i)$ is a Frobenius system for $R\to S$. Hence $S/R$ is a Frobenius
extension.
\end{remarks}

We are now able to give an example of a functor which is full but not naturally full.

\begin{example}\exlabel{3.2}
Let $S$ be a semisimple artinian ring, $T$ a ring, and $0\neq L\in {}_S\Mm_T$.
Consider the ring
$$R=\left(\begin{matrix}
S& L\\ 0& T
\end{matrix}\right),$$
and the obvious projection $\varphi:\ R\to S$ and injection $\nu_R:\ S\to R$.
$\varphi$ is a surjective ring homomorphism, and $\nu_R$ is a left (but not
right) $R$-module section of $\varphi$, so $\varphi_*(S)$ is projective as a left $R$-module.
$S$ is semisimple, so $\varphi_*(N)$ is a projective left $R$-module, for every
$N\in {}_S\Mm$. In particular, the map $\eta_M:\ M\to \varphi_*(\varphi^*(M))$
cosplits in ${}_R\Mm$, so there exists a $\nu_M:\ \varphi_*(\varphi^*(M))\to M$
in ${}_R\Mm$ such that $\eta_M\circ \nu_M= {\varphi_*(\varphi^*(M))}$.
By \prref{2.5}, this means that $\varphi^*$ is full. It follows from \prref{3.1}
that $\varphi^*$ is not naturally full.
\end{example}

\begin{remark}\relabel{3.2}
Consider ring homomorphisms $\alpha:\ A\to B$ and $\beta:\ B\to C$.
If $\alpha_*\circ\beta_*=(\beta\circ\alpha)_*$ is naturally full, then $\beta\circ\alpha$
(and a fortiori $\beta$) is an epimorphism in the category of rings, and
consequently $\beta_*$ is naturally full. This can also be deduced from
\prref{2.3}, since $\alpha_*$ is faithful. Observe that $\alpha_*$ is not
naturally full in general: take the canonical inclusion $\alpha:\ \ZZ\to \ZZ[X]$
and the canonical projection $\beta:\ \ZZ[X]\to \ZZ[X]/(X)$. So the fact that
$HG$ is naturally full does not imply that $H$ is naturally full.
\end{remark}

\begin{remark}\relabel{3.3}
The equivalence $(1)\Leftrightarrow (4)$ in \thref{1.1}
can be obtained from \thref{2.6} as follows. Let
$e=\sum e^1 \otimes\ e^2\in (S\otimes _{R}S)^{S}$. Then, in view
of \equref{3.0.1}, we have:
$$(\xi _{N}^{e}\circ \varepsilon_N)(s\ot_{R}r)=\xi _{N}^{e}(sn)= \sum e^{1}\ot _{R}e^{2}sn.$$
Thus, $ \xi _{N}^{e}$ splits $\varepsilon_N$ iff $$\sum
e^{1}\ot _{R}e^{2}sn=s\ot_{R}n$$ for every $N\in _S\mathcal{M}$,
$n\in N$ and $s\in S$. Clearly this is equivalent to
$1_S\ot_{R}1_S =e\in (S\otimes _{R}S)^{S}$.
\end{remark}

\subsection{Induction and coinduction}\selabel{3.2}
Let $R$ and $S$ be rings and $M\in {}_S\Mm_R$ a $(S,R)$-bimodule.
To $M$ we can associate a pair of adjoint functors
$$F=M\ot_R\bullet:\ {}_R\Mm\to {}_S\Mm~~;~~
G={}_S\Hom(M,\bullet):\ {}_S\Mm\to {}_R\Mm$$
For $Q\in {}_S\Mm$, ${}_S\Hom(M,Q)\in {}_R\Mm$ via the formula
$(m)(r\cdot f)=(mr)f$, for all $m\in M$, $r\in R$ and $f\in {}_S\Hom(M,Q)$. The unit
$\eta$ and the counit $\varepsilon$ of the adjunction are given by
$$\eta_P:\ P\to {}_S\Hom(M,M\ot_R P),~~\eta_P(p)(m)=m\ot_R p;$$
$$\varepsilon_Q:\ M\ot_R {}_S\Hom(M,Q)\to Q,~~\varepsilon_Q(m\ot_R f)=(m)f.$$
Observe that ${}_S\End(M)\in {}_R\Mm_R$ via
$$(m)(r'\cdot f\cdot r)=((mr')f)r;$$
also ${}^*M={}_S\Hom(M,S)\in {}_R\Mm_S$ via
$$(m)(r\cdot f\cdot s)=((mr)f)s.$$

\begin{proposition}\prlabel{3.6}
Let $R$, $S$, $M$, $F$, $G$ be as above. Then
$$\Nat(1_{_S\mathcal{M}}, FG)\simeq (M\ot_R {}^*M)^S.$$
If $M$ is finitely generated and projective as a left $S$-module, then
$$\Nat(GF, 1_{_R\mathcal{M}})\simeq {}_R\Hom{}_R ({}_S\End(M), R).$$
\end{proposition}

\begin{proof}
We refer to \cite{CaenepeelKadison} or \cite{CMZ2002} for full detail.
Let us give a sketch of the proof.
The natural transformation $\xi:\ 1_{_S\mathcal{M}}\to FG$ corresponding to
$\sum_i m_i\ot_R f_i\in (M\ot_R {}^*M)^S$ is the following:
\begin{equation}\eqlabel{3.6.1}
\xi_Q:\ Q\to M\ot_R {}_S\Hom(M,Q),~~(q)\xi_Q=\sum_i m_i \ot_R(?) f_iq
\end{equation}
The natural transformation $\nu:\ GF\rightarrow 1_{_R\mathcal{M}}$
corresponding to an $R$-bimodule map
$E:\ {}_S\End(M)\to R$ is given by
\begin{equation}\eqlabel{3.6.2}
\left( {(\sum_{i} g_i \ot_R
p_i)(\psi_{M,P})}\right)\nu_P=\sum_{i}( g_i)E\cdot p_i
\end{equation}
where
$$\psi_{M,P}:\ {}_S\End(M,)\otimes _R P\stackrel{\sim }{\longrightarrow }{}_S\Hom(M, M\otimes _R
P)$$
$$(m)\left[(g\ot_R p)\psi_{M,P}\right]=(m)g\ot_R p$$ is the canonical
isomorphism of left $R$-modules, as ${}_S M$ is finitely generated
and projective.
\end{proof}

Before we are able to discuss natural fullness of the induction and coinduction
functor, we need a Lemma.

\begin{lemma}\lelabel{3.7}
Let $M\in {}_S\Mm_R$,
$m_1,\ldots,m_n\in M$ and $f_1,\ldots, f_n\in
{}^*M$. Then
\begin{equation}\label{4.15}
m\ot_Rf=\sum_{i=1}^{n} m_i \ot_R(?) f_i\cdot(m)f~~ \text{in}~~
M\ot_R{}_S\Hom (M,Q)
\end{equation}
for every $Q\in {}_S\Mm, f\in {}_S\Hom (M,Q)$ and $m\in
M$ if and only if
\begin{equation}\label{4.14}
m\ot_R {M}=\sum_{i=1}^{n} m_i \ot_R (?)f_i m ~~
\text{in}~~ M\ot_R{}_S\End(M)\,~~\text{for}~\text{every}\quad m\in
M.
\end{equation}
\end{lemma}
\begin{proof}
Taking $Q=M$ and $f={M}$ in (\ref{4.15}), we obtain
(\ref{4.14}). Conversely, assume that (\ref{4.14}) holds. Set
$A={}_S\End(M)$ and consider the canonical isomorphism:
$${\rm can}:\ M\ot_R A\ot_A {}_S\Hom (M,Q)\cong
M\ot_R{}_S\Hom (M,Q)$$
$${\rm can}(m\ot_R a\ot_A
f)=m\ot_Ra\cdot f.$$ By (\ref{4.14}) we have:$$m\ot_R \ot_A
f=\sum_{i=1}^{n} m_i \ot_R(?)f_im\ot_A f.$$ Applying ${\rm can}$,
we find (\ref{4.15}).
\end{proof}

\begin{theorem}\thlabel{3.8}
Let $M$ be an $(S,R)$-bimodule. 1. The coinduction functor
$G={}_S\Hom(M,\bullet)$ is naturally full if and only if there
exists $\sum_{i=1}^{n} m_i \ot_R f_i\in (M\ot_R {}^*M)^S$
satisfying (\ref{4.14}).\\
2. Assume that  $M\in {}_S\Mm$ is finitely generated and
projective. Then the following assertions are equivalent.
\begin{enumerate}
\item The induction functor $F=M\ot_R \bullet$ is
naturally full;
\item there exists $E\in {}_R\Hom{}_R({}_S\End(M), R)$  such that
$\chi\circ E={{}_S\End(M)}$;
\item there exists a central idempotent $e$ of $R$ such that ${}_S\End(M)\cong Re$
 and $\xi:\ R\to {}_S\End(M)\cong Re$ is the projection $\xi(r)=re$.
\end{enumerate}
\end{theorem}
\begin{proof}
1. Let $\xi:\ 1_{_S\mathcal{M}}\rightarrow FG$ be a natural transformation.
By \prref{3.6}, there exists a unique  $\sum_{i=1}^{n} m_i \ot_R f_i\in (M\ot_R {}^*M)^S$
such that $\xi_Q$ is given by \equref{3.6.1}. Then
$$(m\ot_R f)(\xi_Q\circ \varepsilon_Q)=((m)f)\xi_Q=\sum_{i=1}^{n} m_i
\ot_R(?)f_i\cdot(m)f,$$
hence $G$ is a naturally full if and only if (\ref{4.15}) holds,
for every $Q\in {}_S\Mm, f\in
{}_S\Hom (M,Q)$ and $m\in M$. By \leref{3.7}, this is equivalent to (\ref{4.14}).\\
2. Let $\nu:\ GF\rightarrow 1_{_R\mathcal{M}}$ be a natural transformation;
according to \prref{3.6}, there exists a unique  $E\in {}_R\Hom{}_R
({}_S\End(M), R)$
such that $\nu_P$ is given by \equref{3.6.2}. Then we easily compute that
$$\Bigl( {(\sum_{i} g_i \ot_R p_i)(\psi_{M,P})}\Bigr)(\eta_P\circ \nu_P)=
(\sum_{i}( g_i)E\cdot p_i)\eta_P= {(\sum_{i} g_i \ot_R p_i)(\psi_{M,P})}$$
if and only if
$$ (m)(\sum_{i}( g_i)E\cdot p_i)\eta_P=
(m)\Bigl( {(\sum_{i} g_i \ot_R
p_i)(\psi_{M,P})}\Bigr),\text{for}~\text{every}~m\in M$$ if and
only if
$$(m)\ot_R (\sum_{i}( g_i)E\cdot
p_i)=\sum_{i}(m)g_i\ot_R p_i,\text{for}~\text{every}~m\in M$$ if
and only if
$$\sum_{i}m(g_i)E\ot_R p_i=\sum_{i}(m)g_i\ot_R
p_i,\text{for}~\text{every}~m\in M. $$ Therefore $F$ is naturally
full if and only if
$$\eta_P\circ \nu_P={GFP} ~\text{for}~\text{every}~P\in {}_R\Mm$$
if and only if
$$\sum_{i}m(g_i)E\ot_R p_i=\sum_{i}(m)g_i\ot_R
p_i ,$$ for all $ P\in {}_R\Mm$, $p_i\in P$, $g_i \in {}_S\End(M)$
and $m\in M$, if and only if
$$(m)f=mE(f),$$
for all $m\in M$ and $f\in {}_S\End(M)$, if and only if
$$\chi\circ E={{}_S\End(M)}.$$
This proves the equivalence of (1) and (2). The proof of the equivalence of (2) and (3)
is identical to the proof of the equivalence of (2) and (3) in \prref{3.1}.
\end{proof}

\begin{remark}\relabel{3.8.b} \prref{3.1} is a special case of \thref{3.8}: let $i:\ R\to S$
be a ring homomorphism, and view $M=S$ as an $(S,R)$-bimodule.
\end{remark}

\begin{proposition}\prlabel{3.9}
Let $M\in {}_S\Mm_R$ and assume that the
coinduction functor $G={}_S\Hom(M,\bullet)$ is naturally full. Let
$\sum_{i=1}^{n} m_i \ot_R f_i\in (M\ot_R {}^*M)^S$ be as in \thref{3.8}.
Then $e=\sum_{i=1}^{n} (m_i)f_i$ is a central idempotent of
$S$ and
$$sm=sem=esem,$$
for all $s\in S$ and $m\in M$. Therefore $S\simeq S_1\times S_2$ where $S_1=eSe$,
$S$ acts on $M$ via $S_1$ and $M$ is a generator in ${}_{S_1}\Mm$.
\end{proposition}

\begin{proof}
Let $s\in S$.  Since $\sum_{i=1}^{n} m_i \ot_R f_i\in (M\ot_R{}^*M)^S$, we have
$$\sum_{i=1}^{n}s m_i \ot_R f_i=\sum_{i=1}^{n} m_i \ot_R f_i\cdot s.$$
Applying $\varepsilon_S$ to both sides, we find
find
$$\sum_{i=1}^{n}(s m_i )f_i=\sum_{i=1}^{n}( m_i) (f_i\cdot s).$$
Recall that the right $S$-action on ${}^*M$ is given by
$(m)(f\cdot s):=(m)f\cdot s$; using the fact that $f_i$ is left $A$-linear, we find
that
$$\sum_{i=1}^{n}s\cdot(m_i)f_i=\sum_{i=1}^{n}(( m_i)f_i)\cdot s,$$
and this shows that $ e=\sum_{i=1}^{n} (m_i)f_i$ is a central element of $S$.
$e$ is also idempotent since
$$e\cdot e=e\cdot
(1_S)\xi_S\varepsilon_S=(e)\xi_S\varepsilon_S=(1_S)\xi_S\varepsilon_S\xi_S\varepsilon_S=(1_S)\xi_S\varepsilon_S=e,$$
where $\xi_Q$ is given by \equref{3.6.1}.\\
Applying $\varepsilon_M$ to both sides of (\ref{4.14}), we find
$$m=(\sum_{i=1}^{n} (m_i)f_i)m=e\cdot m,$$
hence we have, for all
$s\in S$ and $m\in M$:
$$sm=sem=esem.$$
Taking $Q=S$ in (\ref{4.15}), we find
$$\Tr_M(S):=\sum_{f\in {}_S\Hom(M,S)}\Im(f)=\sum_{f\in {}_S\Hom(M,S)}\sum_ {m\in M}(m)f=\Im(\varepsilon_S).$$
Now from $\xi_S\circ\varepsilon_S={M\ot_R {}^*M}$, it follows that
$$\Im(\varepsilon_S)=\Im(\varepsilon_S\circ\xi_S)=(S)\xi_S\varepsilon_S=Se=eSe.$$
\end{proof}

\begin{proposition}\prlabel{3.10}
Let $M\in {}_S{\mathcal M}_R$. If
$M$ is a generator of ${}_S\Mm$ and $\chi:\ R\rightarrow{}_S\End(M)$ is a ring
epimorphism, then the coinduction functor $G={}_S\Hom(M,\bullet)$
is fully faithful, that is, all counit maps $\varepsilon_Q$ are isomorphisms.
\end{proposition}

\begin{proof}
We will first show that $G$ is naturally full. $M$ is a generator of ${}_S\Mm$,
so there exist $m_1,\ldots,m_n\in M$ and $f_1,\ldots, f_n
\in{^*M}$ such that
$$1_S=\sum_{i=1}^{n}(m_i)f_i.$$
$\chi$ is a ring epimorphism, so it follows from \thref{1.1} that
$$\varepsilon_M :\ M\ot_R {}_S\End(M)\to M, \quad \varepsilon_M
(m\ot_R f) = m\cdot f=(m)f$$
is an isomorphism. From
$$\varepsilon_M (m\ot_R M)=m=1_Sm=\sum_{i=1}^{n}(m_i)f_im=\varepsilon_M
(\sum_{i=1}^{n}m_i\ot_R (?)f_im)),$$
we deduce that
$$m\ot_R M=\sum_{i=1}^{n}m_i\ot_R (?)f_im,$$
for every $m\in M$. Hence
$m_1,\ldots,m_n\in M$ and $f_1,\ldots, f_n$ satisfy
(\ref{4.14}), and (\ref{4.15}) holds, by \leref{3.7}. Now let us prove that
$\sum_{i=1}^{n}m_i\ot_Rf_i\in (M\ot_R {}^*M)^S$, or
\begin{equation}\label{4.18}\sum_{i=1}^{n}sm_i\ot_Rf_i=\sum_{i=1}^{n}m_i\ot_Rf_is,
\end{equation}
for every $s\in S$. We compute
\begin{eqnarray*}
&&\hspace*{-2cm}
\sum_{i=1}^{n} sm_i \ot_R f_i \stackrel{(\ref{4.15}) } {=}
\sum_{i,j=1}^{n}m_j \ot_R (?)f_j\left((sm_i)f_i\right)\\
&=&\sum_{j=1}^{n} m_j \ot_R\sum_{i=1}^{n}(?)f_j\left((sm_i)f_i\right)
=\sum_{j=1}^{n} m_j\ot_R(?)f_js.
\end{eqnarray*}
In fact, for every $m\in M$, we
have:
\begin{eqnarray*}
&&\hspace*{-2cm}
(m)\sum_{i=1}^{n}(?)f_j\left((sm_i)f_i\right)=\sum_{i=1}^{n}(m)f_j\left((sm_i)f_i\right)\\
&=&(m)f_js\sum_{i=1}^{n}  \left((m_i)f_i\right)=(m)f_js=
(m)\left((?)f_js\right),
\end{eqnarray*}
and it follows from \thref{3.8} that $G$ is naturally full, which implies that all the
counit maps $\varepsilon_Q$ split in ${}_S\Mm$.
Since ${}_SM$ is a generator, the
counit map $\varepsilon _Q:\ M\otimes _R{}_S\Hom(M,Q)\to Q$ is an
epimorphism in ${}_S\mathcal{M}$, and therefore it is an isomorphism.
\end{proof}

\subsection{Corings}\selabel{3.3}
Let $R$ be a ring. Recall that an $R$-coring is a comonoid in the monoidal
category ${}_R\Mm_R$. Thus a coring is a triple $(\Cc,\Delta_{\Cc},\varepsilon_{\Cc})$,
where $\Cc$ is an $R$-bimodule, and $\Delta:\ \Cc\to \Cc\ot_R\Cc$ and
$\varepsilon_\Cc:\ \Cc\to R$ are $R$-bimodule maps such that
$$(\Delta_{\mathcal C}\ot_R  {\mathcal C}) \circ \Delta_{\mathcal
C}=({\mathcal C}\ot_R \Delta_{\mathcal C}) \circ \Delta_{\mathcal
C}$$ and
$$(\varepsilon_{\mathcal C} \ot_R
{\mathcal C}) \circ \Delta_{\mathcal C}= ({\mathcal C} \ot_R
\varepsilon_{\mathcal C}) \circ \Delta_{\mathcal C}={\mathcal
C}.$$ $\Delta_\Cc$ is called the comultiplication, and
$\varepsilon_\Cc$ is called the counit. We use the
Sweedler-Heyneman notation
$$\Delta_{\Cc}= c_{(1)}\ot_R c_{(2)},$$
where the summation is implicitely understood. A right $\Cc$-comodule
is a couple $(M,\rho^M)$, where $M$ is a right $R$-module, and
$\rho^M:\ M\to M\ot_R\Cc$ is a right $R$-linear map, called the coaction,
satisfying the conditions
$$(\rho^M \ot_R {\mathcal C}) \circ \rho^M=(M \ot_R
\Delta_{\mathcal C}) \circ \rho^M, \quad (M \ot_R
\varepsilon_{\mathcal C}) \circ \rho^M = M.$$ We use the following
Sweedler-Heyneman notation for coactions:
$$\rho^M(m)= m_{[0]}\ot_Rm_{[1]}$$
Let $M$ and $N$ be right $\Cc$-comodules. A right $R$-linear map
$f:\ M\to N$ is called right $\Cc$-colinear if
$$\rho^N(f(m))=f(m_{[0]})\ot_R m_{[1]},$$
for all $m\in M$. The category of right $\Cc$-comodules and right $\Cc$-colinear
maps is denoted by $\Mm^\Cc_R$.\\
Corings were introduced by Sweedler in \cite{Swe:pre}, and recently revived
by Brzezi\'nski \cite{Brz:cor}. For a detailed study of corings, we refer to
\cite{BrzezinskiWisbauer}.\\
The functor $F:\ \Mm^\Cc_R\to \Mm_R$ forgetting the right
$\Cc$-coaction has a right adjoint $G=\bullet\ot_R\Cc:\Mm_R\to
\Mm^\Cc_R$. For $N\in \Mm_R$, the $\Cc$-comodule structure on
$G(M)=M\ot_R\Cc$ is given by
$$(n\ot_R c)r=n\ot_Rcr~~{\rm and}~~\rho^{G(M)}(n\ot_R c)=n\ot_Rc_{(1)}\ot_Rc_{(2)}$$
The unit and counit of the adjunction are as follows, for $M\in \Mm^\Cc_R$ and
$N\in \Mm_R$:
$$\eta_M=\rho^M:\ M\to M\ot_R\Cc$$
$$\varepsilon_N=I_N\ot_R\varepsilon_\Cc:\ N\ot_R\Cc\to N$$
We will now investigate when $F$ and $G$ are naturally full. In
order to apply \thref{2.6}, we need to know $\Nat(1_{\Mm_R}, FG)$
and $\Nat(GF,1_{\Mm^{\Cc}})$. This computation has been done in
\cite[Proposition 66 and 67]{CMZ2002}. The result is stated in the
next Proposition.

\begin{proposition}\prlabel{3.11}
Let $\mathcal C$ be an $R$-coring, and consider the adjoint pair $(F,G)$
introduced above. Then
$$\Nat(1_{\mathcal{M}_R}, FG)\simeq \Cc^R=\{z\in \Cc~|~rz=zr,~
{\rm for~all~}r\in R\}.$$
The natural transformation $\xi$ corresponding to $z\in \Cc^R$ is the following:
$$\xi_N:\ N\to N\ot_R\Cc,~~\xi_N(n)=n\ot_R z.$$
Also
\begin{eqnarray*}
&&\hspace*{-2cm}
\Nat(GF,1_{\mathcal M_{R}^{\mathcal C}})\simeq
\{\vartheta\in{}_R\Hom_R(\mathcal C\ot_R\mathcal C,R)\\
&|&c_{(1)}\vartheta(c_{(2)}\ot_R d)=\vartheta(c\ot_R
d_{(1)})d_{(2)}{\rm ~for~all~} c,d\in\mathcal C\}
\end{eqnarray*}
The natural transformation $\nu$ corresponding to $\vartheta$ is the following:
$$\nu_M:\ M\ot_R\Cc\to M,~~\nu_M(m\ot_Rc)=m_{[0]}\vartheta(m_{[1]}\ot_R c).$$
\end{proposition}

\begin{proposition}\prlabel{3.12}
Let $\Cc$ be an $R$-coring.
With notation as above, we have:\\
1. The following statements are equivalent:
\begin{enumerate}
\item The functor $G=\bullet\otimes_R \Cc: \Mm_R\to \Mm^{\Cc}_R$
is naturally full; \item there exists  $z\in \Cc^R$ such that
$c=\varepsilon_{\Cc}(c)z$, for all $c\in \Cc$; \item
$\varepsilon_{\Cc}$ splits in $_R\Mm{_R}$, i.e. there is $\xi:R\to
\Cc$ in $_R\Mm{_R}$ such that $\xi\circ \varepsilon_{\Cc}=\Cc$.
\end{enumerate}
2. The following statements are equivalent:
\begin{enumerate}
\item The functor $F:\ \Mm^{\Cc}_R\to \Mm_R$ is naturally full;
\item $c\varepsilon_{\Cc}(d)=\varepsilon_{\Cc}(c)d$, for all
$c,d\in \Cc$; \item $\Delta_\Cc:\Cc\to \Cc\ot_R \Cc$ is surjective
(hence an isomorphism in ${^\Cc_R}\Mm {^\Cc_R}$).
\end{enumerate}
\end{proposition}

\begin{proof}
1. $(1)\Longleftrightarrow (2).$ Take a natural transformation
$\xi:\ 1_{\Mm_R}\to FG$, and the corresponding $z\in \Cc^R$, as in
\prref{3.11}. Then for all $N\in \Mm_R$, $n\in N$ and $c\in \Cc$,
we have that
$$(\xi_N\circ \varepsilon_N)(n\ot_R c)=\xi_N(n\varepsilon_{\Cc}(c))=n\ot_R\varepsilon_{\Cc}(c)z,$$
hence $\xi_N\circ \varepsilon_N={N\ot_R\Cc}$ if and only if
$n\ot_Rc=n\ot_R\varepsilon_{\Cc}(c)z$, for all $n\in N$ and $c\in \Cc$.\\
If $G$ is naturally full, then there exists $\xi$ such that $\xi_N\circ \varepsilon_N={N\ot_R\Cc}$, for all $N$. Taking $N=R$ and $n=1$, we find that
$c=\varepsilon_{\Cc}(c)z$, as needed. The converse is obvious.\\
$(2)\Longrightarrow (3).$  Let $\xi:R\to \Cc$ be defined by
$\xi(r):=rz$.\\
$(3)\Longrightarrow (2).$  Set $z:=\xi(1_R)$.\\
2. $(1)\Longrightarrow (3).$ By \equref{2.6.1} of \thref{2.6},
$\eta_\Cc$ is surjective. Since $\Delta_\Cc=\eta_\Cc$, we
conclude.\\
$(3)\Longrightarrow (2).$ Since $(\varepsilon_\Cc\ot_R \Cc)\circ
\Delta_\Cc=\Cc=(\Cc\ot_R \varepsilon_\Cc)\circ \Delta_\Cc$ this
implies $(\varepsilon_\Cc\ot_R \Cc)=(\Cc\ot_R \varepsilon_\Cc)$
and hence $\varepsilon_\Cc(c)c'=c\varepsilon_\Cc(c')$ for any
$c,c'\in \Cc$.\\
$(2)\Longrightarrow (1).$ Assume that
$c\varepsilon_{\Cc}(d)=\varepsilon_{\Cc}(c)d$, for all $c,d\in
\Cc$, and define $\vartheta:\ \Cc\ot_R\Cc\to R$ by
$$\vartheta(c\ot_Rd)=\varepsilon_{\Cc}(c)\varepsilon_{\Cc}(d).$$
Then
$$c_{(1)}\vartheta(c_{(2)}\ot_R
d)=c \varepsilon_{\Cc}(d)=\varepsilon_{\Cc}(c)d=\vartheta(c\ot_R
d_{(1)})d_{(2)}.$$
The natural transformation $\nu$ corresponding to $\vartheta$ is then a right
inverse of $\eta$, since
$$(\eta_M\circ\nu_M)(m\ot_Rc)=m_{[0]}\varepsilon_{\Cc}(m_{[1]})\varepsilon_{\Cc}(c_{(1)}))\ot_R c_{(2)}=m\ot_Rc,$$
and it follows that $F$ is naturally full.
\end{proof}

\begin{corollary}\colabel{3.12}
Let $\Cc$ be an $R$-bimodule. Then there is a bijection between
the following sets:
\begin{enumerate}
    \item the set of $R$-coring structures on $\Cc$ such that the functor $G=\bullet\otimes_R \Cc: \Mm_R\to
\Mm^{\Cc}_R$ is naturally full;
    \item the set of $R$-ring structures on $\Cc$ such that the restriction of scalars functor is naturally full.
\end{enumerate}
\end{corollary}

\begin{proof}
Let $\Cc$ be an $R$-coring and assume that the functor
$G=\bullet\otimes_R \Cc: \Mm_R\to \Mm^{\Cc}_R$ is naturally full.
By \prref{3.12}, $\varepsilon_{\Cc}$ splits in $_R\Mm{_R}$, i.e.
there is a $\xi:\ R\to \Cc$ in ${}_R\Mm{_R}$ such that $\xi\circ
\varepsilon_{\Cc}=\Cc$. Then we have the following ring structure
on $\Cc$:
$$c\cdot
c'=\xi[\varepsilon_{\Cc}(c)\varepsilon_{\Cc}(c')]~~;~~1_\Cc=\xi(1_R).$$
Moreover, $\xi$ becomes a ring morphism with an $R$-bimodule
section, so the restriction of scalars functor is naturally full, by
\prref{3.1}.\\
Conversely, assume that $\Cc$ is an $R$-ring, such that the
restriction of scalars functor is naturally full. Then there is a
ring morphism $\varphi :\ R\to \Cc$ with a section $E\in
{}_R\Hom_R (\Cc, R)$. We define an $R$-coring structure on $\Cc$
as follows:
$$\Delta(c)=c\ot_R 1_\Cc~~;~~\varepsilon_{\Cc}(c)=E(c).$$
It is straightforward to check that the two constructions are inverse to
each other. For example, take an $R$-coring structure on $\Cc$,
and consider the associated $R$-ring structure. Then
$$c\ot_R 1_\Cc=\sum c_{(1)}\ot_R \varepsilon_{\Cc}(c_{(2)}) 1_\Cc=\sum c_{(1)}\ot_R \xi(\varepsilon_{\Cc}(c_{(2)})) 1_\Cc=\sum c_{(1)}\ot_R
c_{(2)},$$
as needed.
\end{proof}

\begin{corollary}\colabel{3.13}
Let $\Cc$ be an $R$-coring.
1. Assume that $G=\bullet\otimes_R \Cc:\ \Mm_R\to \Mm^{\Cc}_R$
is naturally full. Then we have the following properties.
\begin{enumerate}
    \item The functor $F:\ \Mm^{\Cc}_R\to \Mm_R$ is naturally
    full; the converse property holds if there exists a
$z\in \Cc$ such that $\varepsilon_{\Cc}(z)=1_R$;
    \item $\Cc$ is finitely generated
and projective as a left (right) $R$-module;
    \item the functor $F:\ \Mm {^\Cc_R}\rightarrow
\Mm{_R}$ is a Frobenius functor.
\end{enumerate}
2. If $F:\ \Mm^{\Cc}_R\to \Mm_R$ is naturally full, then it is also
separable.
\end{corollary}

\begin{proof}
1. If $G$ is naturally full, then there exists $z\in \Cc^R$ such
that $c=\varepsilon_{\Cc}(c)z$, for all $c\in \Cc$.\\
(1) We easily compute that
$$c\varepsilon_{\Cc}(d)=\varepsilon_{\Cc}(c)z\varepsilon_{\Cc}(d)=\varepsilon_{\Cc}(c)\varepsilon_{\Cc}(d)z=\varepsilon_{\Cc}(c)d,$$
for all $c,d\in \Cc$, and it follows from \prref{3.12} that $F$ is naturally full.\\
Conversely, assume that $F$ is naturally full, and that there
exists $z\in\Cc$ such that $\varepsilon_{\Cc}(z)=1_R$. Then
$c\varepsilon_{\Cc}(d)=\varepsilon_{\Cc}(c)d$, for all $c,d\in
\Cc$, and, in particular
$c=c\varepsilon_{\Cc}(z)=\varepsilon_{\Cc}(c)z$. $z\in \Cc^R$
since $zr=\varepsilon_{\Cc}(zr)z=\varepsilon_{\Cc}(z)rz=rz$, for
all $r\in R$, and it follows from \prref{3.12} that $G$ is
naturally full.\\
(2) By \prref{3.12}, the functor $G=\bullet\otimes_R \Cc:\ \Mm_R\to
\Mm^{\Cc}_R$ is naturally full if and only if there is $\xi:\ R\to
\Cc$ in $_R\Mm{_R}$ such that $\xi\circ \varepsilon_{\Cc}=\Cc$.
Thus $\Cc$ is finitely generated and projective both as a right
and as a left $R$-module.\\
(3) It is known (see \cite[Theorem 4.1]{Brz:cor}) that the
following statements are equivalent:
\begin{itemize}
\item The forgetful functor $F:\ \Mm {^\Cc_R}\rightarrow \Mm{_R}$ is
Frobenius;
\item $\Cc$ is a finitely generated left
$R$-module and there exists $e\in \Cc^R$ such that the map
$$\phi:\ {}_R\Hom(\Cc,R)\to \Cc,~~ \phi(f)= \sum e_{(1)}f(e_{(2)})$$
 is bijective.
\end{itemize}
Put $e=z$ and observe that
$$\sum z_{(1)}\ot_R z_{(2)}=\sum z_{(1)}\ot_R \varepsilon_{\Cc}(z_{(2)})z=z\ot_R z,$$
so that $\phi(f)=zf(z)$. Consider the map
$$\upsilon:\ \Cc\to
{}_R\Hom(\Cc,R),~~(\upsilon(c))(c')=\varepsilon_{\Cc}(c')\varepsilon_{\Cc}(c).$$
We then compute that
$$(\phi\circ \upsilon)(c)=z((\upsilon)(c))(z)=z\varepsilon_{\Cc}(z)\varepsilon_{\Cc}(c)=z\varepsilon_{\Cc}(z\varepsilon_{\Cc}(c))=z\varepsilon_{\Cc}(c)=c,$$
and
\begin{eqnarray*}
&&\hspace*{-2cm} ((\upsilon\circ
\phi)(f))(c)=((\upsilon(zf(z)))(c)=
\varepsilon_{\Cc}(c)\varepsilon_{\Cc}(zf(z))\\
&=&\varepsilon_{\Cc}(\varepsilon_{\Cc}(c)z)f(z)=\varepsilon_{\Cc}(c)f(z)=f(\varepsilon_{\Cc}(c)z)=f(c).
\end{eqnarray*}
So $\phi$ is bijective and we know from (2) that $\Cc$ is finitely generated and
projective.\\
2. If $F:\ \Mm^{\Cc}\to \Mm_R$ is naturally full,
$\Delta_\Cc:\Cc\to \Cc\ot_R \Cc$ is an isomorphism in ${^\Cc_R}\Mm
{^\Cc_R}$ (the inverse is $\varepsilon_\Cc\ot_R \Cc=\Cc\ot_R
\varepsilon_\Cc$). In particular, $\Cc$ is coseparable and hence,
by \cite[Corollary 3.6]{Brz:cor}, $F$ is separable.
\end{proof}

\begin{examples}\exlabel{3.14}
1. Let $\varphi:\ R \rightarrow S $ be a ring extension. Then ${\Cc}
= S \ot_R S$ is an $S$-coring  with comultiplication $\Delta_{\Cc }:\
S\ot_R S \rightarrow (S \ot_R S) \ot_S (S \ot_R S) $
and counit $\varepsilon_ {\Cc}:\ S \ot_R S \rightarrow S$
given by
$$\Delta_ {\Cc}( a \ot_R a')= (a \ot_R 1_S) \ot_S (1_S
\ot_R a')~~{\rm and}~~\varepsilon_ {\Cc}(a \ot_R a')=aa'.$$
$\Cc$
is called the Sweedler coring associated to the extension
$\varphi:\ R \rightarrow S $. The functor $F$ is naturally full if
and only if the functor $G$ is full if and only if $\varphi$ is a
ring epimorphism.\\
2. Let $I$ be a two-sided idempotent
ideal of a ring $R$ and assume that $I$ is a pure right
$R$-submodule. Then the multiplication map $m_I:I\ot_RI\rightarrow
I$ is bijective and $I$ is an $R$-coring where
$\Delta=m^{-1}:I\rightarrow I\ot_R I$ and
$\varepsilon:I\rightarrow R$ is the inclusion map. It follows from
\prref{3.12} that $F$ is naturally full. By \prref{3.12}, if $G$
is naturally full then $I$ must be finitely generated and
projective as a left (right) $R$-module, but this is not true in
general.
\end{examples}

\begin{remark}
Let $R$ and $S$ be rings and $M\in {}_S\Mm_R$ an
$(S,R)$-bimodule. Then ${}^*M={}_S\Hom(M,S)\in {}_R\Mm_S$. Recall
(see \seref{3.2}) that we have a pair of adjoint
functors
$$F=M\ot_R\bullet:\ {}_R\Mm\to {}_S\Mm~~;~~
G={}_S\Hom(M,\bullet):\ {}_S\Mm\to {}_R\Mm.$$
Now assume that ${}_SM$ is
finitely generated and projective, and let $\{(e_i,{}^*e_i~|~
i=1,\cdots, n\}$ be a finite
dual basis for $M$.
Then we have a coring structure on the $S$-bimodule $\Cc=M\ot_R {}^*M$,
given by
$$\Delta(m\ot_R \mu)= \sum_{i} m\ot_R
{^*}e_i\ot_S e_i\ot_R \mu~~{\rm and}~~
\varepsilon(m\ot_R \mu)=(m)\mu.$$
$\Cc$ is called the comatrix coring associated to $M$. We refer to
\cite{BrzezinskiG,CDV,Kaoutit} for a detailed study of comatrix
corings. We can consider the adjoint pair discussed in \prref{3.11},
in the right handed case:
$$F_{\Cc}:\ {}^\Cc\Mm\to {}_S\Mm~~;~~
G_{\Cc}=\Cc\ot_S\bullet:\ {}_S\Mm\to {}^\Cc\Mm,$$
where $F_{\Cc}$ forgets the left $\Cc$-coaction. Since $M$ is finitely
generated and projective as a left $S$-module, we have for any
left $S$-module that
$${F_{\Cc}}{G_{\Cc}}(N)=M\ot_R {}^*M\ot_S
N\cong M\ot_R {}_S\Hom(M,N)=FG(N).$$
It is easy check that
this isomorphism of left $S$-modules is natural in $N$ and that it
preserves the counits of the two adjunctions. Thus we obtain:
\begin{enumerate}
    \item $\Nat(1_{_S\Mm}, FG)\simeq \Nat(1_{_S\Mm}, {F_{\Cc}}{G_{\Cc}});$
    \item $G$ is naturally full (resp. separable) if and only if $G_{\Cc}$ is naturally full
    (resp. separable).
\end{enumerate}
This explains the similarity between the first parts of Propositions \ref{pr:3.6} and
\ref{pr:3.11}. Moreover it tells that $G$ is naturally full if and only if
there exists $z\in (M\ot_R {}^*M)^S$ such that $m\ot_R \mu=z\cdot
(m)\mu$, for all $m\ot_R \mu\in M\ot_R {}^*M$ (apply the left sided version
of \prref{3.12}).\\
Now let us look at the $R$-bimodule $A={}^*M\ot_S M\cong {}_S\End(M)$.
Thus $A$ is a ring, with multiplication and counit given by
$$(\mu\ot_S m)\cdot
(\tau\ot_S t)=\mu\ot_S (m)\tau\cdot t~~;~~1_A=\sum_i{^*e_i}\ot_S e_i,$$
and we have a ring morphism
$$\varphi:\ R\to A,~~\varphi(r)= 1_A\cdot r= r\cdot 1_A.$$
We then have a pair of adjoint functors
$$F_A=A\ot_R\bullet:\ {}_R\Mm\to {}_A\Mm~~;~~G_A:\ {}_A\Mm\to {}_R\Mm,$$
where $G_A$ is restriction of scalars. Since $M$ is finitely generated and
projective as a left $S$-module, we have, for any $N\in {}_R\Mm$,
$${G_A}{F_A}(N)={}^*M\ot_S M\ot_R
N\cong {}_S\Hom(M,M\ot_R N)=GF(N).$$
This isomorphism in natural in $N$ and preserves the units of the adjunctions.
Hence we obtain the following:
\begin{enumerate}
    \item $\Nat(GF, 1_{{}_R\Mm})\simeq \Nat(F_AG_A, 1_{{}_R\Mm});$
    \item $F$ is naturally full (resp. separable) if and only if $F_A$ is naturally full
    (resp. separable).
\end{enumerate}
This explains the similarity between the second parts of Propositions \ref{pr:3.6} and
\ref{pr:3.11}. Moreover it clarifies the analogy between the second
parts of \prref{3.1} and \thref{3.8}.
\end{remark}

Let $\Cc$ be an $R$-coring. $g\in\Cc$ is called grouplike if $\Delta_\Cc(g)
=g\ot_R g$ and $\varepsilon_\Cc(g)=1$. There is a bijective correspondence
between the grouplike elements of $\Cc$
and right (or left) $\Cc$-comodule structures on $R$,
see \cite{Brz:cor}: the right $\Cc$-coaction on $R$ corresponding to $g$ is
given by $\rho(r)=1\ot_Rgr$.\\
Fix a grouplike element $g\in \Cc$, and define the coinvariants $M^{{\rm co}\Cc}$
of $M$ by
$$M^{{\rm co}{\Cc}}=\{m \in M \mid \rho_M (m)=m \ot_R g \}.$$
In particular,
$$B=R^{{\rm co}{\Cc}}=\{r \in R~|~rg=g r \}$$
is a subring of $R$ and $M^{{\rm co}{\Cc}}\in \Mm_B$. We obtain a functor
$\tilde{G}=(\bullet)^{{\rm co}\Cc}:\ \Mm^\Cc_R\to \Mm_B$, which has a left adjoint
$\tilde{F}=\bullet\ot_BR$. For any $N\in \Mm_B$, $N\ot_BR$ is a right $\Cc$-comodule
via the right $\Cc$-coaction on $R$. The unit and counit of the adjunction are
the following, for every $N\in \Mm_B$ and $M\in \Mm^\Cc_R$:
$$\alpha_N:\ N\to (N\ot_BR)^{{\rm co}\Cc},~~\alpha_N(n)=n\ot_B 1;$$
$$\beta_M:\ M^{{\rm co}\Cc}\ot_B R\to M,~~\beta_M(m\ot_B a)=ma.$$

\begin{proposition}\prlabel{3.15}
Let $\Cc$ be an $R$-coring with a fixed grouplike element $g$.
If there exists an $(R,R)$-bimodule map $\chi:\ \Cc\ot_R\Cc\to R$
such that
\begin{equation}\eqlabel{3.15.1}
c_{(1)}\chi(c_{(2)}\ot_R d)=\chi(c\ot_R d_{(1)})d_{(2)},
\end{equation}
for all $c,d\in \Cc$ and
$$\chi(g\ot_R g)=1,$$
then the functor $\tilde{F}=\bullet\ot_BR:\ \Mm_B\to\Mm^\Cc_R$
is fully faithful, which means that it is separable and naturally full.
\end{proposition}

\begin{proof}
We first show that the map
$$t:\ R\to B,~~t(r)=\chi(gr\ot_R g)$$
is a $B$-bimodule map. $t$ is left $B$-linear since
$$t(br)=\chi(gbr\ot_R g)=\chi(bgr\ot_R g)=b\chi(gr\ot_R g)=bt(r),$$
for all $b\in B$ and $r\in R$. A similar argument shows that $t$ is right $B$-linear.\\
Let us show that $t(r)\in B$, for all $r\in R$: applying \equref{3.15.1} with
$c=gr$, $d=g$, we find that
$g\chi(gr\ot_R g)=\chi(gr\ot_R g)g$, hence $t(r)=\chi(gr\ot_R g)\in B$.
Also observe that $t^2=t$, since $t(b)=b$, for $b\in B$.\\
Take $N\in \Mm_B$. We claim that the map
$$\theta_N:\ (N\ot_B R)^{{\rm co}\Cc}\to N,~~
\theta_N(\sum_i n_i\ot_B
r_i)=\sum_in_it(r_i)=\sum_in_i\chi(g\ot_Rr_ig)$$ is inverse to
$\alpha_N$. Obviously $(\theta_N\circ \alpha_N)(n)=n$, for all
$n\in N$.\\Take $\sum_i n_i\ot_B r_i\in (N\ot_B R)^{{\rm co}\Cc}$.
Then
$$\sum_i n_i\ot_B gr_i=\sum_i n_i\ot_B r_ig,$$
so
$$\sum_i n_i\ot_B gr_i\ot_R g=\sum_i n_i\ot_B r_ig\ot_Rg,$$
and
$$\sum_i n_i\ot_B \chi(gr_i\ot_R g)=\sum_i n_i\ot_B \chi(r_ig\ot_Rg).$$
Since $\chi(gr_i\ot_R g)\in B$ and $\chi$ is left $R$-linear, we find
$$\sum_i n_i \chi(gr_i\ot_R g)\ot_B 1_R=\sum_i n_i\ot_B r_i\chi(g\ot_Rg),$$
and, since $\chi(g\ot_Rg)=1$,
$$(\alpha_N\circ\theta_N)(\sum_i n_i\ot_B r_i)=\sum_i n_i\ot_B r_i.$$
\end{proof}

Recall from \cite{Brz:cor} that a coring $\Cc$ is called coseparable if
the forgetful functor $F:\ \Mm^\Cc_R\to \Mm^R$ is separable; this is equivalent
to the existence of an $(R,R)$-bimodule map
$\xi:\ \Cc\ot_R\Cc\to R$ satisfying \equref{3.15.1}
such that $\chi\circ \Delta_\Cc=\varepsilon_\Cc$. If $\Cc$ is coseparable, then
the conditions of \prref{3.15} are satisfied, and $\tilde{F}$ is fully faithful, for
every choice of the grouplike element $g$.

\subsection{Homomorphisms of corings}\selabel{3.4}
Let $\mathcal C$ be an
$R$-coring and let $\mathcal D$ be an $S$-coring. A coring
homomorphism from $\Cc  $ to $\Dd  $ is a pair
$(\Phi,\varphi)$, where $\varphi:\ R\rightarrow S$ is a ring
homomorphism and $\Phi:\ \Cc  \rightarrow \Dd  $ is a
homomorphism of $R$-modules such that:
$$\omega_\mathcal{D,D}\circ
(\Phi\ot_R \Phi)\circ \Delta_\Cc  =\Delta_\Dd  \circ
\Phi~~;~~\varphi\circ
\varepsilon_\Cc  =\varepsilon_\Dd  \circ \Phi,$$ where
$\omega_\mathcal{D,D}:\Dd  \ot_R \Dd  \rightarrow
\Dd  \ot_S \Dd  $ is the canonical map induced by
$\varphi$. \\
We have a functor $S\ot_R
\bullet:{^\Cc  _R}\Mm \rightarrow {^\Dd  _S}\Mm$,
where the left $\Dd  $-comodule structure of $S\ot_R M$ is
given by the map $s\ot_R m \mapsto \sum s \Phi(m_{[-1]})\ot_S 1_S
\ot_R m_{[0]}$. Observe that $S\ot_R \Cc  $ is a right
$\Cc  $-comodule via $S\ot_R \Delta_\Cc  $.
In a similar way, we have a functor $\bullet\ot_R
S:\ \Mm{^\Cc  _R} \rightarrow \Mm{^\Dd  _S}$.\\
Let $L$ be an $(S,R)$-bimodule, and let $\rho_L^\Dd  :\
L\rightarrow L\ot_S \Dd  $ and ${^\Cc  }\rho_L:\ L\rightarrow \Cc
\ot_R L$ be comodule structures on $L$. It is clear that ${^\Cc
}\rho_L$ is a morphism of right $\Dd  $-comodules if and only if
$\rho_L^\Dd  $ is a morphism of left $\Cc  $-comodules. In this
case we say that $L$ is a $(\Cc  ,\Dd  )$-bicomodule; we can then
define the category ${}^\Cc\Mm ^\Dd$ of $(\Cc  ,\Dd
)$-bicomodules. The morphisms are $(R,S)$-bimodule maps that are
left $\Cc  $-colinear and right $\Dd $-colinear.\newline

Let $f,g:X\to Y$ be a pair of morphisms of right modules over a
ring $R$ and let $i:E\to X$ be its equalizer (that is the kernel
of $f-g$). We will say that a left $R$-module $_RM$ preserves the
equalizer of $(f,g)$ if the map $i\ot_R M:E\ot_R M\to X\ot_R M$ is
the equalizer of the pair $(f\ot_R M,g\ot_R M)$. Analogously for
the right-hand side.

\begin{remark}\relabel{3.19}
\prref{3.20} and the subsequentTheorems are inspired by \cite{Tore}. The
referee pointed out to us that it is not sufficient, as in the original
statement in \cite{Tore}, to assume that $_R\Cc$  preserves equalizers. In fact,
one needs the assumption that $_R(\Cc\ot_R\Cc) $ preserves equalizers.
More precisely, let $M=N\ot _S S \ot _R \Cc$ and let $(K,k)$
be the equalizer of $\rho_N^{\Dd}\ot _S S\ot _R\Cc$ and $N\ot
_S{^{\Dd}\rho _{S\ot _R \Cc}}$. One has to prove that $(K,k)$
inherits a $\Cc  $-bicomodule structure. Now, by means of the
universal property of the equalizers, there is a unique morphism
$\rho_K^{\Cc}:K\to K\ot_R\Cc$ such that
\begin{equation*}
(k\ot_R \Cc)\rho_K^{\Cc}=\rho_M^{\Cc}k.
\end{equation*}
Using the coassociativity of the coaction of $M$ and the relation
above, we obtain:
\begin{eqnarray*}
&&\hspace*{-15mm}
(k\otimes _{R}\Cc\otimes _{R}\Cc)(\rho ^{\Cc}_{K}\otimes
_{R}\Cc)\rho^{\Cc}_{K}= [(k\otimes _{R}\Cc)\rho^{\Cc} _{K}\otimes _{R}\Cc]\rho^{\Cc} _{K} \\
&=&(\rho^{\Cc} _{M}k\otimes _{R}\Cc)\rho^{\Cc} _{K} =
(\rho^{\Cc} _{M}\otimes _{R}\Cc)(k\otimes _{R}\Cc)\rho^{\Cc} _{K}
\end{eqnarray*}
\begin{eqnarray*}
&=&(\rho^{\Cc} _{M}\otimes _{R}\Cc)\rho ^{\Cc}_{M}k =
(M\otimes _{R}\Delta _{\Cc})\rho ^{\Cc}_{M}k \\
&=&(M\otimes _{R}\Delta _{\Cc})(k\otimes _{R}\Cc)\rho ^{\Cc}_{K} =
 (k\otimes _{R}\Cc\otimes _{R}\Cc)(K\otimes _{R}\Delta
_{\Cc})\rho ^{\Cc}_{K}.
\end{eqnarray*}
Thus, since $_R(\Cc\ot_R\Cc) $ preserves equalizers, the morphism
$k\otimes _{R}\Cc\otimes _{R}\Cc$ is a monomorphism so that we get
\begin{equation*} (\rho ^{\Cc}_{K}\otimes
_{R}\Cc)\rho^{\Cc}_{K}=(K\otimes _{R}\Delta _{\Cc})\rho
^{\Cc}_{K}.
\end{equation*}
Note also that $_R\Cc$ preserves equalizers, whenever
$_R(\Cc\ot_R\Cc)$ does. Therefore this requirement is a little
stronger and it has replaced the first one in all the following
results.
\end{remark}

\begin{proposition}\prlabel{3.20} (cf. \cite[Proposition 5.4]{Tore})
The $(S,R)$-bimodule $S\ot_R \Cc  $ is a $(\Dd  ,\Cc
)$-bicomodule. Moreover, if $_R(\Cc\ot_R\Cc)  $ preserves the
equalizer of $(\rho_N^\Dd  \ot_S S\ot_R \Cc  ,N\ot_S
{}^\Dd\rho_{S\ot_R \Cc })$ for every right $\Dd  $-comodule $N$,
then the functor
$$G=\bullet\square_{\Dd  }(S\ot_R \Cc  ):\Mm {^\Dd  _S}
\rightarrow \Mm {^\Cc  _R}$$ is right adjoint to
$$F=\bullet\ot_R S:\Mm {^\Cc  _R}\rightarrow \Mm
{^\Dd  _S}.$$
\end{proposition}

The unit $\eta:\ 1_{\Mm {^\Cc  _R}}\to GF$ and the counit
$\varepsilon:\ FG\to 1_{\Mm {^\Dd  _S}}$ are given by
$$\eta_M(m)=\sum (m_{[0]}\ot_R
1_S)\square_{\Dd  }(1_S \ot_R m_{[1]})$$
and
$$\varepsilon_N((n\square_{\Dd  }(s\ot_R c))\ot_R
s'):=ns\varphi(\varepsilon_\Cc  (c))s'$$
for all $M\in {\Mm
{^\Cc  _R}}$ and  $N\in {\Mm{^\Dd  _S}}$.

\begin{theorem}\thlabel{3.20}
Assume that $_R(\Cc\ot_R \Cc)$ preserves the equalizer of
$(\rho_N^\Dd \ot_S S\ot_R \Cc  ,N\ot_S {}^\Dd \rho_{S\ot_R \Cc })$
for every $N\in \Mm^\Dd  _S$, and that $M_R$ preserves the
equalizer of $(\rho_{\Cc  \ot_R S}^\Dd  \ot_S S\ot_R \Cc  ,{\Cc
\ot_R S}\ot_S {}^\Dd\rho_{S\ot_R \Cc  })$ for every $M\in \Mm^\Cc
_R$. The functor $F=\bullet\ot_R S:\ \Mm^\Cc _R\rightarrow
\Mm{^\Dd  _S}$ is naturally full if and only if $\eta_\Cc  :\ \Cc
\to GF\Cc  $ cosplits in ${}^\Cc  _R\Mm {^\Cc  _R}$, i.e. there is
a homomorphism of $\Cc  $-bicomodules $\nu_\Cc  :\ GF\Cc  \to \Cc
$ such that $\eta_\Cc  \circ \nu_\Cc  ={GF\Cc  }$.
\end{theorem}

\begin{proof}
By \thref{2.6}, $F$ is naturally full if and only if $\eta:\
1_{\Mm {^\Cc  _R}}\to GF$ cosplits, i.e. there exists a natural
transformation $\nu :\ GF \to 1_{\Mm {^\Cc  _R}}$ such that
$\eta_M\circ \nu_M = {GFM}$ for all $M \in {\Mm {^\Cc  _R}}$. In
this case, obviously $\eta_\Cc  \circ \nu_\Cc  ={GF\Cc  }$ and
$\nu_\Cc  $ is the
required map.\\
Conversely, assume that there is a homomorphism of $\Cc
$-bicomodules $\nu_\Cc  :\ GF\Cc  \to \Cc  $ such that $\eta_\Cc
\circ \nu_\Cc  ={GF\Cc  }$. Let $r_M:\ M\ot_R R\to M$ be the
canonical isomorphism. According to \cite[Theorem 5.6]{Tore}, we
have a natural transformation $\nu :\ GF \to 1_{\Mm {^\Cc _R}}$
defined as follows:
$$\nu_M=r_M\circ (M\ot_R
\varepsilon_{\Cc}\circ\nu_{\Cc})\circ \kappa_M,$$ where
$\kappa_M:\ GFM\to M\ot_R GF\Cc$ is a natural isomorphism (see the
proof of \cite[Theorem 5.6]{Tore}). It is also proved in
\cite{Tore} that that $\kappa_M\circ \eta_M=(M\ot_R
\eta_{\Cc})\circ \rho_M^\Cc  $. In a similar way, we can show that
$$(\rho_M^\Cc  \ot_R GF\Cc)\circ \kappa_M=(M\ot_R
\kappa_\Cc)\circ \kappa_M.$$
Hence we have:
\begin{eqnarray*}
&&\hspace*{-15mm} \rho_M^\Cc  \circ \nu_M=\rho_M^\Cc  \circ r_M\circ
(M\ot_R \varepsilon_{\Cc}\circ \nu_{\Cc})\circ \kappa_M\\
&=&(M\ot_R r_\Cc)\circ (\rho_M^\Cc  \ot_R R)\circ (M\ot_R \varepsilon_{\Cc}\circ \nu_{\Cc})\circ \kappa_M\\
&=&(M\ot_R r_\Cc)\circ (M\ot_R \Cc \ot_R \varepsilon_{\Cc}\circ \nu_{\Cc})\circ (\rho_M^\Cc  \ot_R GF\Cc)\circ \kappa_M\\
&=&(M\ot_R r_\Cc)\circ (M\ot_R \Cc \ot_R \varepsilon_{\Cc}\circ
\nu_{\Cc})\circ (M\ot_R \kappa_\Cc)\circ \kappa_M\\
&=&(M\ot_R \nu_{\Cc})\circ \kappa_M.
\end{eqnarray*}
Now let us compute:
\begin{eqnarray*}
&&\hspace*{-2cm}\kappa_M\circ \eta_M\circ \nu_M=(M\ot_R \eta_{\Cc})\circ
\rho_M^\Cc  \circ \nu_M\\&=&(M\ot_R \eta_{\Cc})\circ (M\ot_R
\nu_{\Cc})\circ \kappa_M=\kappa_M.
\end{eqnarray*}
$\kappa_M$ is a monomorphism, so it follows that $\eta_M\circ \nu_M =
{GFM}$.
\end{proof}

\begin{theorem}\thlabel{3.21}
Assume that $_RS$ and $_R\Cc  $ preserve the
equalizer of $(\rho_N^\Dd  \ot_S S\ot_R \Cc  ,N\ot_S
{}^\Dd\rho_{S\ot_R \Cc  })$ for every $N\in
\Mm^\Dd  _S$. The functor
$G=\bullet\square_{\Dd  }(S\ot_R \Cc  ):\Mm
{^\Dd  _S} \rightarrow \Mm {^\Cc  _R}$ is naturally
full if and only if the $\Dd$-bicomodule map
$$\hat{\Phi}:\ S\ot_R
\Cc\ot_R S\to \Dd,~~\hat{\Phi}(s\ot_R c\ot_R s)= s\Phi(c)s'$$
splits in
$^\Dd  _S\Mm {^\Dd  _S}$,that is, there is a homomorphism
of $\Dd  $-bicomodules $\hat{\Psi}:\ \Dd\to S\ot_R \Cc\ot_R S$
such that $\hat{\Psi}\circ \hat{\Phi}={S\ot_R \Cc\ot_R S}$.
\end{theorem}

\begin{proof}
By \thref{2.6}, $G$ is naturally full if and only if
$\varepsilon:\ FG\to 1_{\Mm^{\Dd}_S}$ splits, that is, there
exists a natural transformation $\xi:\ 1_{\Mm^{\Dd}_S}\to FG$ such
that $\xi_N\circ \varepsilon_N = {FGN}$ for all $N \in
{\Mm^{\Dd}_S}$. In particular $\xi_\Dd\circ \varepsilon_\Dd =
{FG\Dd}$. In \cite[5.6]{Tore}, it is proved that
$$\hat{\Phi}=\varepsilon_\Dd\circ ({^\Dd}\rho_{S\ot_R \Cc}\ot_R S).$$
Let $l_N:\ S\ot_S N\to N$ be the canonical isomorphism, for every
$N\in {\Mm^{\Dd}_S}$. Then we can write
\begin{eqnarray*}
&&\hspace*{-15mm}(l_{S\ot_R \Cc }\ot_R S)\circ (\varepsilon_\Dd\ot_S(S\ot_R
\Cc)\ot_R
S)\circ \xi_\Dd\circ \hat{\Phi}\\
&=& (l_{S\ot_R \Cc }\ot_R S)\circ (\varepsilon_\Dd\ot_S(S\ot_R
\Cc)\ot_R S)\circ ({^\Dd}\rho_{S\ot_R \Cc}\ot_R S)\\
&=& (l_{S\ot_R \Cc }\ot_R S)\circ (l_{S\ot_R \Cc}^{-1}\ot_R
S)={S\ot_R \Cc \ot_R S},
\end{eqnarray*}
so we can choose $\hat{\Psi}=(l_{S\ot_R \Cc }\ot_R S)\circ (\varepsilon_\Dd\ot_S(S\ot_R
\Cc)\ot_R S)\circ \xi_\Dd$.\\
Conversely, assume that there is a homomorphism of $\Dd
$-bicomodules $\hat{\Psi}:\Dd\to S\ot_R \Cc\ot_R S$ such that
$\hat{\Psi}\circ \hat{\Phi}={S\ot_R \Cc\ot_R S}$. In \cite[Theorem
5.8]{Tore}, it is proved that the map
$$N\rTo^{\rho_N^{\Dd}} N\ot_S \Dd\rTo{N\ot_S \hat{\Psi}}N\ot_S S\ot_R \Cc\ot_R S$$
factorizes through a natural transformation
$$\xi_N:\ N\to FGN=(N\sq_\Dd (S\ot_R \Cc))\ot_R S.$$
Now let $i:GN\to N\ot_S (S\ot_R \Cc)$, for every $N\in \Mm^\Dd
_S$, be the equalizer of $(\rho_N^\Dd  \ot_S S\ot_R \Cc ,N\ot_S
{^\Dd  }\rho_{S\ot_R \Cc })$ which is preserved by $_RS$. Then the
morphism $i\ot_R S$ is a monomorphism. In \cite{Tore} we also find
that: $$(i\ot_R S)\circ \xi_N=(N\ot_S \hat{\Psi})\circ
\rho_N^{\Dd}.$$ It is easy to check (just apply the definition of
the cotensor product) that:
$$\rho_N^{\Dd}\circ
\varepsilon_N=(N\ot_S \hat{\Phi})\circ (i\ot_R S).$$ Thus, we
have:
\begin{eqnarray*}
&&\hspace*{-2cm} (i\ot_R S)\circ \xi_N\circ \varepsilon_N=(N\ot_S
\hat{\Psi})\circ
\rho_N^{\Dd}\circ \varepsilon_N\\
&=& (N\ot_S \hat{\Psi})\circ (N\ot_S \hat{\Phi})\circ (i\ot_R
S)=i\ot_R S.
\end{eqnarray*}
Since $i\ot_R S$ is a monomorphism, we conclude that
$\xi_N\circ \varepsilon_N={FGN}.$
\end{proof}

\begin{examples}
Consider the adjunction
$(F=\bullet\ot_R S,G=\bullet\square_{\Dd  }(S\ot_R\Cc  ))$ between the
categories $\Mm {^\Cc  _R}$ and $\Mm{^\Dd  _S}$.
Observe first that the category
$\Mm^{\Cc}=\Mm_R$ in the case where $\Cc=R$.\\
1) Let  $\Cc$ and $\Dd$ be coalgebras over a field $\Kk$ and
$(\Phi,\varphi)=(\Phi,\mathbb{K})$. Then  $F$ is the corestriction of
coscalars functor, and $G=\bullet\square_{\Dd}\Cc$.
By \thref{3.20}, $F$ is naturally full if and only if
$\eta_\Cc:\Cc\to GF\Cc=\Cc\sq_\Dd \Cc$, $\eta_\Cc(c)= \sum c_{(1)}\sq_\Dd
c_{(2)}$ cosplits as a $\Cc$-bicomodule map. Since
$$(\Phi\sq_\Dd \Cc)\circ \eta_\Cc=\Cc=(\Cc\sq_\Dd \Phi)\circ \eta_\Cc,$$
we have that $\Phi\sq_\Dd \Cc=\Cc\sq_\Dd \Phi$ and hence
$\varepsilon_\Cc(c)c'=c\varepsilon_\Cc(c')$ for any $c,c'\in \Cc$
(recall that we have a canonical isomorphism $\Dd\sq_\Dd
\Cc\to\Cc,~~d\sq_\Dd c\mapsto \varepsilon_\Dd(d)c$).
In view of \prref{3.12}, this last property means that the forgetful
functor $\Mm^{\Cc}_\Kk\to \Mm_\Kk$ is naturally full. Now, if
$\Cc\neq 0$, there is a $c\in \Cc$ such that $\varepsilon_\Cc(c)\neq
0$. Since $\mathbb{K}$ is a field, we can put
$z:=\varepsilon_\Cc(c)^{-1}c$ so that $\varepsilon_\Cc(z)=1$. By
\coref{3.13}, the functor $G=\bullet\otimes_\Kk \Cc: \Mm_\Kk\to
\Mm^{\Cc}_{\Kk}$ is naturally full and, by \coref{3.12},
$\Cc$ has a $\Kk$-ring structure by means of a ring homomorphism
$\xi:\Kk\to \Cc$ such that $\xi\circ \varepsilon_\Cc=\Cc$. Since
$\Kk$ is a field, $\xi$ is an isomorphism.\\
By \thref{3.21}, $G$ is naturally full if and only if the $\Dd$-bicomodule map
$$\hat{\Phi}:\ \Kk\ot_\Kk \Cc\ot_\Kk \Kk\to \Dd,~~\hat{\Phi}(k\ot_\Kk c\ot_\Kk
k')= k\Phi(c)k'$$
splits as a $\Dd$-bicomodule map. Since
$\Kk\ot_\Kk \Cc\ot_\Kk \Kk\simeq \Cc$, this is equivalent to
$\Phi$ splitting as a $\Dd$-bicomodule map.\\

2) Let $(\Phi,\varphi)=(\varphi,\varphi)$, that is, $\Cc=R$ and
$\Dd=S$. In this case the adjunction $(F,G)$ reduces to extension and
restriction of scalars:
$$F=\bullet\ot_R
S:\ \Mm {_R}\rightarrow \Mm {_S}~~;~~G:\ \Mm {_S}
\rightarrow \Mm {_R}.$$
By \thref{3.20}, $F$ is naturally full if and only if $\eta_R:\ R\to
GFR=R\ot_R S$, $\eta_R(r)= r\ot_R 1_S$, cosplits in ${}_R\Mm{_R}$. Since
$l_S:\ R\ot_R S\to S$, $l_S(r\ot_R s)= \varphi(r)s$ is an isomorphism
and $l_s\circ \eta_R=\varphi$, this is equivalent to
$\varphi$ cosplitting in ${}_R\Mm{_R}$, that is, there exists an $E\in
{}_R\Hom (S, R){}_R$ such that $ \varphi\circ E={S}$. We
have therefore recovered condition 2-2) of \prref{3.1}. \\
By \thref{3.21}, $G$ is naturally full if and only if the
$S$-bimodule map $\hat{\varphi}:\ S\ot_R R\ot_R S\to S,~~\hat{\varphi}(s\ot_R r\ot_R
s)= s\varphi(r)s'$, splits in ${}_S\Mm {_S}$. Since
$\hat{\varphi}\circ (S\ot_R l_s^{-1})=\varepsilon_S$, the counit of
the Sweedler coring $S\ot_RS$, this is equivalent to
$\varepsilon_S$ splittin in ${}_S\Mm {_S}$. Since $\tau:\ S\to S\ot_R
S$, $\tau(s)= 1_S\ot_R s$ is a retraction of $\varepsilon_S$ (i.e.
$\varepsilon_S\circ \tau=S$), we can conclude that $G$ is
naturally full if and only if $\varepsilon_S$ is injective. So we have
recovered condition (5) of \thref{1.1}.\\
3) Let $(\Phi,\varphi)=(\varepsilon_{\Cc},R)$, that is
$\Dd=S=R$. The adjoint pair $(F,G)$ is now the adjoint pair from
\prref{3.12}. \\
By \thref{3.20}, $F$ is naturally full if and only if
$\eta_\Cc=\Delta_\Cc:\ \Cc\to GF\Cc=\Cc\ot_R \Cc$ cosplits in
$^\Cc  _R\Mm {^\Cc  _R}$. Since $(\varepsilon_\Cc\ot_R
\Cc)\circ \Delta_\Cc=\Cc=(\Cc\ot_R \varepsilon_\Cc)\circ
\Delta_\Cc$ this is equivalent to $\Delta_\Cc$ being
surjective, and we recover condition (3) of
\prref{3.12}-2.\\
By \thref{3.21}, $G$ is naturally full if and only if the
$R$-bimodule map $\hat{\varepsilon}_{\Cc}:\ R\ot_R \Cc\ot_R R\to
R$, $\hat{\varepsilon}_{\Cc}(r\ot_R c\ot_R r')= r\varepsilon_{\Cc}(c)r'$, splits in
${}_R\Mm{_R}$. Since $R\ot_R \Cc\ot_R R\cong \Cc$, this is
equivalent to  $\varepsilon_{\Cc}$ splitting in ${}_R\Mm{_R}$: we
recover condition (3) of \prref{3.12}-1.
\end{examples}

\begin{remark}
Let $(R,\Cc)_\psi$ be an entwining structure over a commutative
ring $K$. It was shown in \cite[Proposition 2.2]{Brz:cor} that
$R\ot_K \Cc$ has a coring structure in such a way that the
category of comodules over $R\ot_K \Cc$ is isomorphic to the
category of entwined modules $\Mm^{\Cc}_R(\psi)$. Every morphism
$(\varphi,\Theta):\ (R,\Cc)_\psi\to (S,\Dd)_\phi$ of entwining
structures induces a coring homomorphism
$(\varphi\ot_K\Theta,\varphi):\ R\ot_K \Cc\to S\ot_K \Dd$.
Applying Theorems \ref{th:3.20} and \ref{th:3.20}, we can give
necessary and sufficient conditions for the faithful fullness of
the functors in the adjoint pair $(F=\bullet\ot_R S, G=
\bullet\square_{\Dd}\Cc$ between $\Mm^{\Cc}_R(\psi)$ and
$\Mm^{\Dd}_R(\phi)$.
\end{remark}

\begin{center}
{\sc Acknowledgement}
\end{center}
We thank Jos\'{e} (Pepe) G\'{o}mez Torrecillas  for
providing \exref{3.14} 2), and the referee for his comments on \prref{3.20}.

\end{document}